\documentclass[final]{siamltex}

\usepackage{amsmath}
\usepackage{graphicx}
\usepackage{latexsym}
\usepackage{amssymb}
\usepackage{amsbsy}
\usepackage{dsfont}
\DeclareSymbolFont{msbm}{U}{msb}{m}{n}
\DeclareMathSymbol{\C}{\mathalpha}{msbm}{'103}
\DeclareMathSymbol{\R}{\mathalpha}{msbm}{'122}
\DeclareMathSymbol{\Z}{\mathalpha}{msbm}{'132}
\DeclareMathSymbol{\N}{\mathalpha}{msbm}{'116}

\newcommand{\D}{\mathcal{D}}

\newcommand{\F}{\mathcal{F}}

\newtheorem{remark}{Remark}

\def\F{{\cal F}}

\def\RR{\mathbb R}
\def\e{\varepsilon}

\def\nuu{s}
\def\mum{\cal M}

\def\be{\begin{equation}}
\def\ee{\end{equation}}
\def\bea{\begin{eqnarray}}
\def\ba{\begin{array}{l}\displaystyle}
\def\eea{\end{eqnarray}}
\def\ea{\end{array}}

\def\proof{{\bf Proof \par}}

\begin{document}
\title{Implicit-Explicit linear multistep methods for stiff kinetic equations}

\author{Giacomo Dimarco\thanks{Mathematics and Computer Science Department, University of
Ferrara, Ferrara, Italy ({\tt giacomo.dimarco@unife.it}).}
\and Lorenzo Pareschi\thanks{Mathematics and Computer Science Department, University of
Ferrara, Ferrara, Italy ({\tt lorenzo.pareschi@unife.it}).}
}
\maketitle

\begin{abstract}
We consider the development of high order asymptotic-preserving linear multistep methods for kinetic equations and related problems. The methods are first developed for BGK-like kinetic models and then extended to the case of the full Boltzmann equation. The behavior of the schemes in the Navier-Stokes regime is also studied and compatibility conditions derived. We show that, compared to IMEX Runge-Kutta methods, the IMEX multistep schemes have several advantages due to the absence of coupling conditions and to the greater computational efficiency. The latter is of paramount importance when dealing with the time discretization of multidimensional kinetic equations.
\end{abstract}

\maketitle



\section{Introduction}
The development of robust numerical methods for the solution of kinetic equations with stiff collision terms has been a very active field of research in the recent years \cite{BLM, BJ, CJR, degond1, dimarco2, dimarco5, Filbet3, Jin2, Lem, KT}. Implicit-Explicit (IMEX) Runge-Kutta methods and exponential schemes have provided effective computational tools to solve such equations in stiff regimes \cite{boscarino,dimarco6,dimarco7,Filbet,HLP,PRimex,Puppo,toscani,Qin}. For a comprensive introduction to the subject we refer the reader to the recent survey \cite{Acta}.

This paper is the natural continuation of this research direction on effective numerical methods for stiff kinetic equations and considers the development of IMEX multistep methods for Boltzmann-type equations. IMEX multistep methods were originally proposed in \cite{Ascher2} and subsequently further developed and studied in \cite{DB, FHV, HR, HRS}. Comparison between IMEX multistep methods and IMEX Runge-Kutta methods were presented in \cite{HR}. Let us also mention that fully implicit multistep methods for BGK equations were considered in \cite{GRS}.  


We consider kinetic equations of the form~\cite{cercignani}
\be
\partial_t f + v\cdot\nabla_{x}f
=\frac1{\varepsilon}Q(f), \label{eq:1b}\ee
where $\varepsilon$ is the Knudsen number a non dimensional quantity directly proportional to the mean free pat between particles.
Here $f(x,v,t)$ is a non negative function
describing the time evolution of the distribution of particles with
velocity $v \in \R^{d_v}$, $d_v \geq 1$ and position $x \in \Omega \subset \R^{d_x}$, $d_x \geq 1$
at time $ t
> 0$. 

For notation simplicity in the sequel we will omit the
dependence of $f$ from the independent variables $x,v,t$ unless
strictly necessary. The operator $Q(f)$ characterizes the particles
interactions, in the case of the quadratic Boltzmann collision operator of rarefied gas dynamics $Q(f)=Q_B(f,f)$ where for $d_v\geq 2$ we have
 \be Q_{B}(f,f)=\int_{\RR^{d_v}\times S^{d_v-1}} B(|v-v_*|,n)
[f(v')f(v'_*)-f(v)f(v_*)]\,dv_*\,dn \label{eq:Q} \ee where \be
v'=v+\frac12(v-v_*)+\frac12|v-v_*|n,\quad
v'_*=v+\frac12(v-v_*)-\frac12|v-v_*|n, \ee and $B(|v-v_*|,n)$ is a
nonnegative collision kernel characterizing the details of the
collision. 

The operator $Q(f)$ is such that the local
conservation properties are satisfied \be\int_{\R^{d_v}} \phi(v) Q(f)\,
dv=:\langle \phi Q(f)\rangle=0 \label{eq:QC}\ee where
$\phi(v)=\left(1,v,{|v|^2}/{2}\right)^T$ are the collision invariants. In
addition it satisfies the entropy inequality \be
\frac{d}{dt}H(f) = \int_{\R^{d_v}} Q(f)\log f dv
\leq 0,\qquad H(f)=\int_{\R^{d_v}}f\log f\,dv. \label{eq:entropy} \ee
The functions such that
$Q(f)=0$ are the local Maxwellian equilibrium functions
\be M[f]=M(\rho,u,T)=\frac{\rho}{(2\pi
T)^{d_v/2}}\exp\left(\frac{-|u-v|^{2}}{2T}\right), \label{eq:M}\ee
where $\rho$, $u$, $T$ are the density, mean velocity and
temperature of the gas in the $x$-position and at time $t$ defined as
\be (\rho,\rho u,E)^T=\langle \phi f \rangle, \qquad
T=\frac1{d_v\rho}(E-\rho|u|^2). \ee Due to its computational complexity, the Boltzmann collision operator
$Q_{B}(f,f)$ is often replaced in applications by simpler operators, like the BGK operator which substitutes the binary interactions with a
relaxation towards the equilibrium of the form~\cite{BGK} \be\label{BGK}
Q_{BGK}(f)=\mu(M[f]-f),\ee
where $\mu > 0$ is usually assumed to be proportional to the macroscopic density $\rho(x,t)$. The validity of this operator in
describing the physics of non equilibrium phenomena has been the
subject of many papers in the past~\cite{BGK,BPESBGK,cercignani}.

In this paper, we first consider the development of IMEX multistep methods for BGK-like kinetic models, namely models for which the implicit solution of the stiff collision term can be solved explicitly and therefore will not increase the overall computational cost of the scheme. Next we extend the schemes to the case of the full Boltzmann equation by adopting a penalization strategy analogous to the one introduced in \cite{dimarco7, dimarco8, Filbet, Qin} for Runge-Kutta and exponential methods. In particular, we analyze the behavior of the schemes in the Navier-Stokes regime and show that the methods are capable to described correctly this asymptotic behavior. We emphasize that a similar analysis for IMEX Runge-Kutta and exponential methods is still lacking. We refer to \cite{BLM, Filbet} for some consistency results with the Navier-Stokes regime for some simple one step methods. In addition, compared to IMEX Runge-Kutta methods, the IMEX multistep schemes have several advantages due to the absence of coupling conditions and to the greater computational efficiency. The latter is of paramount importance when dealing with the time discretization of multidimensional kinetic equations. On the other hand a particular care should be used in the choice of the penalization factor for high order methods due to the stronger stability constraints in comparison with IMEX Runge-Kutta methods and exponential methods. 

The rest of the paper is organized as follows. First in Section 2 we
recall some basic aspects on kinetic equations and their fluid dynamic limits.
The notion of asymptotic preservation is also introduced. In
Section 3, we consider IMEX linear multistep
schemes applied to kinetic equations and derive conditions for  asymptotic preservation and
asymptotic accuracy in Navier-Stokes regimes. Next, in section 4 we
introduce the penalized IMEX linear multistep schemes for the full Boltzmann model. We analyze their asymptotic preservation and the behavior in the Navier-Stokes regime. Monotonicity in the homogeneous case and stability of the penalized approach are also discussed. Finally, in Section 5, several numerical examples confirms our theoretical findings. Some concluding remarks are reported in Sections 6.

\section{Hydrodynamic limits and asymptotic-preserving methods}

For notation simplicity in this Section we consider only the case $d_v=3$. Integrating (\ref{eq:1b}) against the collision
invariants in the velocity space leads to the following set of non
closed conservations laws \be
\partial_t \langle \phi f\rangle+{\rm div}_x
\langle v\otimes\phi f\rangle=0.\label{eq:macr}\ee

Close to fluid regimes, the mean free path between two collisions is
very small. In this situation, passing to the limit $\varepsilon\rightarrow 0$ we formally obtain $Q(f)=0$ from (\ref{eq:1b}) and so $f=M[f]$. Thus, at least formally, we
recover the closed hyperbolic system of compressible Euler equations
\be
\partial_t U+{\rm div}_x \F(U)=0
\label{eq:Euler} \ee with
$U=\langle \phi M[f]\rangle = (\rho,\rho u,E)^T$ and
\[
\F(U)=\langle v\otimes\phi M[f]\rangle=\left(
\begin{array}{c}
\rho u  \\
\rho u \otimes u+pI  \\
Eu+pu     
\end{array}
\right)
,\quad p=\rho T,\quad T=\frac{1}{3}\left(\frac{2E}{\rho}-|u|^2\right),
\]
where $I$ is the identity matrix. Note that the above conclusions are
independent on the particular choice of $Q(f)$ provided it satisfies (\ref{eq:QC}) and admits Maxwellian of the form (\ref{eq:M}) as local equilibrium functions.

For small but non zero values of the Knudsen number, the evolution
equation for the moments can be derived by the so-called
Chapman-Enskog expansion~\cite{BGP00, cercignani}.
This originates the compressible Navier-Stokes equations
as a second order approximation with respect to $\varepsilon$ to the
solution of the Boltzmann equation
\be
\partial_t U+{\rm div}_x \F(U)=\varepsilon\, {\rm div}_x \D(\nabla_x U)
\label{eq:NavierStokes} \ee with
\be
\D(\nabla_x U)=
\left(
\begin{array}{c}
  0 \\
\nu\sigma(u)\\
\kappa\nabla_x T+\nu\sigma(u)\cdot u     
\end{array}
\right),
\quad 
\sigma(u)=\frac12\left(\nabla_x u +(\nabla_x u)^T -\frac2{3}{\rm div}_x uI\right),
\label{eq:NSsigma}
\ee
and the viscosity $\nu$ and the thermal conductivity $\kappa$ are defined according to the
linearized Boltzmann operator with respect to the local Maxwellian \cite{cercignani, Golse}. The Prandtl number is the ratio $Pr = 5\nu/(2\kappa)$.
Note that, the
particular choice of the collision operator $Q(f)$ influences the structure
of the limiting Navier-Stokes system in terms of the corresponding Prandtl number \cite{BGP00, BPESBGK, Caflisch, cercignani, Golse}.

%

The construction of numerical schemes which are capable to capture the Euler limit just described is closely connected with the notion of asymptotic-preserving schemes. Here, in agreement with~\cite{Jin2, PRimex} we give the following definition of
asymptotic preserving methods for equation (\ref{eq:1b})
\begin{definition}
A consistent time discretization method for (\ref{eq:1b}) of stepsize
$\Delta t$ is {\em asymptotic preserving (AP)} if, independently of
the initial data and of the stepsize $\Delta t$, in
the limit $\varepsilon\to 0$ becomes a consistent time
discretization method for the reduced system (\ref{eq:Euler}).
\end{definition}

This definition does not imply that the scheme preserves the order
of accuracy in time in the stiff limit $\varepsilon\to 0$. In the
latter case we will say that the scheme is \emph{asymptotically accurate (AA)}. One of the main question which is left open by the above definition is the behavior of the numerical method in the Navier-Stokes regime.  



\section{IMEX linear multistep schemes for kinetic equations}

First we introduce the general formulation of IMEX linear multistep schemes for kinetic equations together with some preliminary definitions. 


For a kinetic equation of the form (\ref{eq:1b}) we consider schemes based on a combination of $\nuu$-step explicit and implicit linear multistep methods
\begin{equation}
   f^{n+1} = - \sum_{j=0}^{\nuu-1} a_j f^{n-j} - \Delta t \sum_{j=0}^{\nuu-1} b_j\, v\cdot\nabla_x f^{n-j}+ \Delta t \sum_{j=-1}^{\nuu-1} c_j \frac{1}{\varepsilon}Q(f^{n-j}),
\label{eq:GIMEX1}
\end{equation}
where $c_{-1} \neq 0$. Methods for which $c_j=0$, $j=0,\ldots,\nuu-1$ are referred to as implicit-explicit backward differentiation formula, IMEX-BDF in short. Another important class is represented by implicit-explicit Adams methods, for which $a_0=-1$, $a_j=0$, $j=1,\ldots,\nuu-1$.

We refer to~\cite{Ascher2, HR, FHV, HairerWanner} for more details on
the order conditions for IMEX multistep schemes. Let us recall that an order $p$ scheme is obtained provided that
\begin{equation}
\begin{aligned}
&1+\sum_{j=0}^{\nuu-1} a_j =0,\\
1-\sum_{j=1}^{\nuu-1}& j a_j =\sum_{j=0}^{\nuu-1} b_j =\sum_{j=-1}^{\nuu-1} c_j,\\
\frac12+\sum_{j=1}^{\nuu-1} \frac{j^2}{2}a_j&=-\sum_{j=1}^{\nuu-1} jb_j=c_{-1}-\sum_{j=1}^{\nuu-1}j c_j\\
\vdots\\
\frac1{p!}+\sum_{j=1}^{\nuu-1} \frac{(-j)^p}{p!} a_j=-\sum_{j=1}^{\nuu-1}&\frac{(-j)^{p-1}}{(p-1)!} jb_j=\frac{c_{-1}}{(p-1)!}+\sum_{j=1}^{\nuu-1}\frac{(-j)^{p-1}}{(p-1)!} c_j.
\end{aligned}
\label{eq:GIMEXcond}
\end{equation}
Moreover the following theorem holds true \cite{Ascher2}
\begin{theorem}
For the $\nuu$-step IMEX scheme (\ref{eq:GIMEX1}) we have
\begin{enumerate}
\item If $p\leq \nuu$ the $2p+1$ constraints of (\ref{eq:GIMEXcond}) are linearly independent, therefore there exist $\nuu$-step IMEX multistep schemes of order $\nuu$.
\item A $\nuu$-step IMEX multistep scheme has accuracy at most $\nuu$.
\item The family of $\nuu$-step IMEX multistep schemes of order $\nuu$ has $\nuu$ parameters. 
\end{enumerate}
\end{theorem}

In the sequel 
we will make use of the following vector notation 
\begin{equation}
   f^{n+1} = - a^T \cdot F - \Delta t\, b^T \cdot L(F) + \frac{\Delta t}{\varepsilon} c^T \cdot Q(F)+\frac{\Delta t}{\varepsilon} c_{-1} Q(f^{n+1}),
\label{eq:GIMEXv}
\end{equation}
where $a=(a_0,\ldots,a_{\nuu-1})^T$, $b=(b_0,\ldots,b_{\nuu-1})^T$, $c=(c_{0},\ldots,c_{\nuu-1})^T$, $F=(f^n,\ldots,f^{n-\nuu+1})^T$, $L(F)=(v\cdot \nabla_x f^n,\ldots,v\cdot \nabla_x f^{n-\nuu+1})^T$,  and $Q(F)=(Q(f^{n}),\ldots,Q(f^{n-\nuu+1}))^T$ are $\nuu$-dimensional vectors. In Table \ref{tb:examples} we report some examples of IMEX multistep schemes \cite{Ascher2,HR}. 


\begin{remark}
An interesting feature of IMEX-BDF schemes is that they can be rewritten in splitting form as
\bea
\label{eq:IMEXBDFs1}
   f^{n+1/2} &=& - a^T \cdot F - \Delta t\, b^T \cdot L(F) \\
   f^{n+1} &=& f^{n+1/2} + \frac{\Delta t}{\varepsilon} c_{-1} Q(f^{n+1}),
\label{eq:IMEXBDFs2}
\eea
where the collision step (\ref{eq:IMEXBDFs2}) corresponds to a simple application of the backward Euler scheme, whereas the multistep approach is now limited to the transport step (\ref{eq:IMEXBDFs1}). \end{remark}

\subsection{Asymptotic preserving IMEX multistep schemes}

In this section we give conditions for an IMEX multistep scheme to satisfy asymptotic
preservation and asymptotic accuracy. Here we do not consider the computational challenges related to the inversion of the implicit collision operator $Q(f)$. We will focus on these aspects in Section 4.

We can state the following theorem which show that BDF-IMEX schemes are asymptotic preserving and
asymptotically accurate. 
In order to do this we first introduce the notion of vector of initial
steps {consistent} with the limit problem.
\begin{definition} The vector of initial steps $F$ in the IMEX multistep scheme 
(\ref{eq:GIMEXv}) is said
\emph{consistent} or \emph{well-prepared} if
\be f^{n-j}(x,v)=M[f^{n-j}(x,v)]+g^{n-j}_{\varepsilon}(x,v),\qquad \lim_{\varepsilon\to 0} g^{n-j}_{\varepsilon}(x,v)=0,\quad j=0,\ldots\nuu-1.\label{eq:cid}\ee
\end{definition}
The above definition is satisfied, for example, using an asymptotic-preserving IMEX Runge-Kutta scheme of the type developed in \cite{dimarco7} as initializing method. 

\begin{table}[t]
\caption{Examples of multistep IMEX schemes up to order $5$}
\begin{center}
{\small
\begin{tabular}{l|c c c}
\hline\\[-.25cm]
Scheme & $a^T$ & $b^T$ & $(c_{-1},c^T)$\\
\hline\\[-.25cm]
IMEX-BDF1 & $-1$ & $1$ & $\left(1,0\right)$\\[+.25cm]
IMEX-CN2 & $\left(-1,0\right)$ & $\left(\frac32,-\frac12\right)$ & $\left(\frac1{2},\frac12,0\right)$\\[+.25cm]
IMEX-MCN2 & $\left(-1,0\right)$ & $\left(\frac32,-\frac12\right)$ & $\left(\frac9{16},\frac38,\frac1{16}\right)$\\[+.25cm]
IMEX-BDF2 & $\left(-\frac43,\frac13\right)$ & $\left(\frac43,-\frac23\right)$ & $\left(\frac23,0,0\right)$\\[+.25cm]
IMEX-SG2 & $\left(-\frac34,0,-\frac14\right)$ & 
$\left(\frac32,0,0\right)$ & $\left(1,0,0,\frac12\right)$\\[+.25cm]
IMEX-BDF3 & $\left(-\frac{18}{11},\frac{9}{11},-\frac{2}{11}\right)$ & $\left(\frac{18}{11},-\frac{18}{11},\frac{6}{11}\right)$ & $\left(\frac{6}{11},0,0,0\right)$\\[+.25cm]
IMEX-AD3 & $\left(-1,0,0\right)$ & $\left(\frac{23}{12},-\frac{4}{3},\frac{5}{12}\right)$ & $\left(\frac{4661}{10000},\frac{15551}{30000}, \frac{1949}{30000},-\frac{1483}{30000}\right)$\\[+.25cm]
IMEX-TVB3 & {$\left(-\frac{3909}{2048},\frac{1367}{1024},-\frac{873}{2048}\right)\!\!$} & {$\left(\frac{18463}{12288},-\frac{1271}{768},\frac{8233}{12288}\right)\!\!$} & {$\left(\frac{1699}{12288},\frac{1089}{2048},-\frac{1139}{12288},-\frac{367}{6144}\right)$}\\[+.25cm]
IMEX-BDF4 & $\left(-\frac{48}{25},\frac{36}{25},-\frac{16}{25},\frac3{25}\right)$ & $\left(\frac{48}{25},-\frac{72}{25},\frac{48}{25},-\frac{12}{25}\right)$ & $\left(\frac{12}{25},0,0,0,0\right)$\\[+.25cm]
IMEX-TVB4 & {$\left(-\frac{21531}{8192},\frac{22573}{8192},\right.$}& {$\left(\frac{13261}{8192},-\frac{75029}{24576},\right.$}& 
{$\left(\frac{4207}{8192}, 
-\frac{3567}{8192},
\frac{697}{24576},\right.$}\\[+.25cm]
& {$\left.-\frac{12245}{8192},\frac{2831}{8192}\right)\!\!$} & {$\left.\frac{54799}{24576},-\frac{15245}{24576}\right)\!\!$}  & {$\left.\frac{4315}{24576},
-\frac{41}{384}\right)$}
\\[+.25cm]
IMEX-BDF5 & $\left(-\frac{300}{137},\frac{300}{137},-\frac{200}{137},\frac{75}{137},-\frac{12}{137}\right)$ & $\left(\frac{300}{137},-\frac{600}{137},\frac{600}{137},-\frac{300}{137},\frac{60}{137}\right)$ & $\left(\frac{60}{137},0,0,0,0,0\right)$\\[+.25cm]
IMEX-TVB5 & $\left(
-\frac{13553}{
4096},
\frac{38121}{
8192},
-\frac{7315}{
2048},\right.$
&
$\left(
\frac{10306951}{
5898240},
-\frac{13656497}{
2949120},
\frac{1249949}{
245760},
\right.\!\!\!\!$
&
$\left(
\frac{4007}{
8192},
-\frac{4118249}{
5898240},
\frac{768703}{
2949120},
\right.$
\\[+.25cm]

& $\left.\frac{6161}{
4096},
-\frac{2269}{
8192}
\right)$
 & 
$\left.
-\frac{7937687}{
2949120},
\frac{3387361}{
5898240}
\right)$
&
$\left.
\frac{47849}{
245760},
-\frac{725087}{
2949120},
\frac{502321}{
5898240}
\right)$
\\[+.25cm]
\hline
\end{tabular}
}
\end{center}
\label{tb:examples}
\end{table}

\begin{theorem}
\label{th:ap} If the vector of initial steps is well-prepared, in the limit
$\varepsilon\rightarrow 0$, scheme
(\ref{eq:GIMEXv}) becomes the explicit multistep scheme characterized by
($a, b$) applied to the limit
Euler system (\ref{eq:Euler}).
\end{theorem}

 \proof To prove the above theorem let us first multiply the IMEX
multistep scheme (\ref{eq:GIMEXv}) by the collision
invariant vector $\phi(v)=(1, v, \frac12{|v|^2})^T$ and integrate the result in velocity
space. We obtain the explicit multistep method applied to the
moment system (\ref{eq:macr})
\begin{equation}
\label{eq:GIMEXfm} \\
\displaystyle   \langle \phi f^{n+1}\rangle = -\sum_{j=0}^{s-1} a_j\,\langle\phi f^{n-j} \rangle - \Delta t\, \sum_{j=0}^{s-1} b_j\, {\rm div}_x \langle v\otimes \phi f^{n-j}\rangle.
\end{equation}
Now let us rewrite equation (\ref{eq:GIMEXv})
in the following form 
\begin{equation}
   \varepsilon f^{n+1} = - \varepsilon a^T \cdot F - \varepsilon \Delta t\, b^T \cdot L(F) + \Delta t\, c^T \cdot Q(F) +\Delta t\, c_{-1} Q(f^{n+1}),
   \label{eq:espiIMEX}
\end{equation}
which, as $\varepsilon\to 0$ reduces to
\begin{equation}
   0 = c^T \cdot Q(F) + c_{-1} Q(f^{n+1}).
   \label{eq:espiIMEX2}
\end{equation}
Since the initial steps are well-prepared, as $\varepsilon\rightarrow 0$ we get $f^{n-j}=M[f^{n-j}]$, $j=0,\ldots,\nuu-1$ which implies $Q(F)\equiv 0$ and therefore since $c_{-1}\neq 0$ we have
\be
Q(f^{n+1})=0 \Rightarrow f^{n+1}=M[f^{n+1}].
\ee
Thus (\ref{eq:GIMEXfm})
 becomes the explicit multistep method applied to the limiting Euler system (\ref{eq:Euler})
\be
\label{eq:GIMEXfeuler} 
U^{n+1} = -\sum_{j=0}^{s-1} a_j\, {U}^{n-j} - \Delta t\, \sum_{j=0}^{s-1}b_j\, {\rm div}_x \F(U)^{n-j},
\ee
where 
\be
\F(U)^{n-j}=\left(
\begin{array}{c}
\rho^{n-j} u^{n-j}  \\
\rho^{n-j} u^{n-j}\otimes u^{n-j}+p^{n-j}I \\
E^{n-j} u^{n-j}+p^{n-j} u^{n-j}    
\end{array}
\right).
\ee
$\Box$

Note that for an IMEX-BDF scheme we have $c\equiv 0$ and therefore, even for non well-prepared initial conditions, in
the limit $\varepsilon\rightarrow 0$ equation (\ref{eq:espiIMEX2}) reduces to $Q(f^{n+1})=0$ and  the distribution function is
projected over the equilibrium $f^{n+1}\rightarrow M[f^{n+1}]$. Thus, after $\nuu$ time steps, the IMEX-BDF scheme becomes equivalent to the explicit multistep scheme (\ref{eq:GIMEXfeuler}). 

Thus, we can state
\begin{theorem}
For arbitrary initial steps an IMEX-BDF scheme, in the limit
$\varepsilon\rightarrow 0$, after $\nuu$ time steps becomes the explicit multistep scheme characterized by
($a, b$) applied to the limit
Euler system (\ref{eq:Euler}).
\end{theorem}

\begin{remark}
For IMEX-BDF schemes, in case the initial steps are not well-prepared the numerical solution may exhibit an initial layer which gives rise to a reduction of accuracy in the numerical solution. This phenomena can be cured using smaller time steps or extrapolation techniques only for the very first steps of the computation.
\end{remark}

Finally, let us consider the property of the schemes to preserve a stationary state of the system. In the field of PDEs this is usually referred as \emph{well-balanced property} and involves not only the time discretization of the system but also the space discretization (see \cite{Gosse} for a recent survey on the topic). 
More precisely, a desirable property of multistep IMEX schemes is that they preserve a steady state solution such that \be\varepsilon v\cdot \nabla_x f = Q(f),\ee for a fixed value of $\varepsilon$. It is easy to show the following result.

\begin{theorem}
If the vector of initial steps is well-balanced, namely $F=f^n e$ with $\varepsilon v\cdot \nabla_x f^n = Q(f^n)$, then the multistep IMEX scheme (\ref{eq:GIMEXv})
preserves the steady state solution or equivalently $f^{n+1}=f^n$.
\end{theorem}

\proof
It is enough to observe that under the assumptions we have 
\[
\varepsilon L(F)=Q(F)
\] 
or $\varepsilon (v\cdot \nabla_x f^n) e=Q(f^n)e$ with $e=(1,\ldots,1)^T$. We can write
\[
f^{n+1}-\frac{\Delta t}{\varepsilon} c_{-1} Q(f^{n+1}) = -a^T\cdot f^n e - \frac{\Delta t}{\varepsilon}( b^T-c^T)\cdot Q(f^n)e,
\]
since from the order conditions $a^T\cdot e=-1$, $c_{-1}+c^T\cdot e=b^T\cdot e$ the conclusion follows.\\
$\Box$

\subsection{IMEX linear multistep schemes for relaxation operators}
\label{par:BGK}
In this paragraph we consider the particular case of BGK relaxation operators of the form $Q_{BGK}(f)=\mu(M[f]-f)$, where $\mu=\mu(x,t)$ depends only on the macroscopic quantities.
A fundamental property of the IMEX scheme
(\ref{eq:GIMEXv}) applied to relaxation operators is that it can be solved explicitly
despite the nonlinearity of $M[f]$.

In this case (\ref{eq:GIMEXv}) takes the form
\begin{equation}
   f^{n+1} = - a^T \cdot F - \Delta t\, b^T \cdot L(F) + \frac{\Delta t}{\varepsilon} c^T \cdot {\mum}(M[F]-F)+\frac{\Delta t}{\varepsilon} c_{-1} \mu^{n+1}(M[f^{n+1}]-f^{n+1}),
\label{eq:GIMEXbgk}    
\end{equation}   
where ${\mum}={\rm diag}\{\mu^n,\ldots,\mu^{n-\nuu+1}\}$ and $M[F]=(M[f^n],\ldots,M[f^{n-s+1}])^T$, which can be rewritten as   
\bea
\nonumber
f^{n+1}   =&& \frac{\varepsilon}{{\varepsilon+c_{-1}\mu^{n+1}\Delta t}}\left( - a^T \cdot F - \Delta t\, b^T \cdot L(F)\right) \\[-.25cm]
\label{eq:GIMEXbgk2}   
\\[-.25cm]
\nonumber
&+& \frac{\Delta t}{{\varepsilon+c_{-1}\mu^{n+1}\Delta t}}\left({c^T} \cdot {\mum}(M[F]-F)+ c_{-1}\mu^{n+1}M[f^{n+1}]\right)
   \eea
where the only implicit term is 
$\mu^{n+1}M[f^{n+1}]$ which depends only on the
moments $\langle\phi f^{n+1}\rangle$. If we now integrate equation
(\ref{eq:GIMEXbgk}) against the collision invariants thanks to the
conservations (\ref{eq:QC}) we obtain the explicit moment scheme (\ref{eq:GIMEXfm}).
Thus $\langle \phi
f^{n+1}\rangle$, and so $\mu^{n+1}M[f^{n+1}]$, can be explicitly evaluated
and the numerical solution (\ref{eq:GIMEXbgk2}) can be directly computed. 


\subsubsection{The space-homogeneous case}
Monotonicity properties for IMEX linear multistep schemes have been studied in~\cite{FHV, HR, HRS}. Here, we restrict our analysis to the
space homogeneous case ($L(f)=0$) by focusing on two properties which are of particular interest for kinetic equations, namely non negativity of the solution and entropy inequality.
As we will see, even in this simplified setting, the resulting conditions for implicit linear multistep schemes are rather restrictive. 


In the homogeneous case the method reduces to a simple application of the implicit scheme and reads
\begin{eqnarray}
\nonumber
   f^{n+1} &=& -\lambda a^T \cdot F  + (1-\lambda)\left(\frac{c^T}{c_{-1}} \cdot (M[F]-F)+ M[f^{n+1}]\right)\\[-.25cm]
     \label{eq:s1}
   \\[-.25cm]
   \nonumber
   &=& -\left(\lambda a^T+(1-\lambda) \frac{c^T}{c_{-1}}\right)\cdot F+(1-\lambda)\left(\frac{c^T}{c_{-1}} \cdot e + 1\right) M[f^n],
\end{eqnarray}   
where we have set $e=(1,\ldots,1)^T$, $\lambda={\varepsilon}/({\varepsilon+c_{-1}\mu\Delta t})\in [0,1)$ and we used the fact that in the homogeneous case the Maxwellian $M[f]$ and $\mu$ are  independent of time.
If we now define $z=\mu\Delta t/\varepsilon$ we have $\lambda=1/(1+c_{-1}z)$ and we can state the following
\begin{proposition}
Sufficient conditions to guarantee that $f^{n+1}\geq 0$ when $F\geq 0$ in (\ref{eq:s1}) are that
\begin{eqnarray}
\label{eq:posn}
a^T +z c^T \leq 0,\\
c^T\cdot e+ c_{-1} \geq 0.
\label{eq:posn1}
\end{eqnarray}
\end{proposition}

Note that conditions (\ref{eq:posn}) must be
interpreted component by component. In particular, (\ref{eq:posn})
depend on $z$ and originates the time step restriction. 
For BDF schemes $c^T\equiv 0$ and conditions reduce to $a_j \leq 0$, $\forall\,j$ which are never satisfied except for the first order backward Euler method which is the unique method unconditionally monotone. For Adams Moulton methods conditions (\ref{eq:posn}) become $z \leq 1/c_0$, $c_0 > 0$ and $c_j \leq 0$, $j=1,\ldots,\nuu-1$.
Examples of methods satisfying these restrictions are the popular Cranck-Nicolson scheme $a^T=(-1,0)$, $c^T=(1/2,0)$, $c_{-1}=1/2$ for $z \leq 2$, 
and the third order scheme $a^T=(-1,0,0)^T$, $c^T=(2/3,-1/12,0)^T$, $c_{-1}=5/12$ for $z\leq 3/2$. Implicit Adams schemes of order higher then third never satisfy $c_j \leq 0$, $j=1,\ldots,\nuu-1$. Another example of second order nonnegative method is the implicit part of IMEX-SG2 scheme $a^T=(-3/4,0,-1/4)$, $c^T=(0,0,1/2)$, $c_{-1}=1$ for $z \leq 1/2$.


Let us remark that since from the first order conditions $a^T\cdot e = -1$ we can also write
\be
f^{n+1}=-\frac{a^T+zc^T}{1+c_{-1}z}\cdot F+\left(\frac{e^T}{\nuu}+\frac{a^T+zc^T}{1+c_{-1}z}\right)\cdot eM.\label{eq:simpp}\ee
Therefore we can state
\begin{proposition}
If conditions (\ref{eq:posn})-(\ref{eq:posn1}) are satisfied 
then (\ref{eq:simpp}) is a convex combination of the initial steps and the Maxwellian state.
\end{proposition}

The above proposition gives the following entropy inequality (see (\ref{eq:entropy})) 
\be
H(f^{n+1}) \leq \sum_{j=0}^{\nuu-1} \alpha_j H(f^{n-j})+\left(1-\sum_{j=0}^{\nuu-1}\alpha_j\right) H(M) \leq \max_{0\leq j \leq s-1} H(f^{n-j}),
\ee
where 
\[
\alpha_j = -\frac{a_j+zc_j}{1+c_{-1}z}.
\]

\begin{remark}
Its is clear that there are few implicit linear multistep methods satisfying conditions (\ref{eq:posn})-(\ref{eq:posn1}). More relaxed conditions can be derived for linear multistep methods if we take into account the starting procedure that generates the vector of $f^{n-j}$ \cite{HRS}. Here we do not explore further this direction.
\end{remark}

\subsubsection{Navier-Stokes asymptotics}
A fundamental problem when dealing with AP schemes for kinetic equations is the behavior of the method in the Navier-Stokes regime. In order to do this we rewrite the general multistep IMEX scheme of order $p$ as
\begin{equation}
   \frac{f^{n+1}+ a^T \cdot F}{\Delta t} +  b^T \cdot L(F) = \frac{1}{\varepsilon} c^T \cdot {\mum}(M[F]-F)+\frac{1}{\varepsilon} c_{-1} \mu^{n+1}(M[f^{n+1}]-f^{n+1}),
   \label{eq:IMEX_BGK}
\end{equation}
and consider the discrete Chapman-Enskog expansion taking
\be
f^{n+1}=M[f^{n+1}]+\varepsilon g^{n+1},\quad F=M[F]+\varepsilon G,
\ee
where  $\langle \phi g^{n+1}\rangle=0$ and $\langle \phi G\rangle \equiv 0$ with $G=(g^n,\ldots,g^{n-s+1})^T$. 

Inserting the above expansions in the numerical method yields
\bea
\nonumber
   \frac{M[f^{n+1}]+ a^T \cdot M[F]}{\Delta t} &+&  b^T \cdot L(M[F]) \\[-.25cm]
   \label{eq:IM_CE}
   \\[-.25cm]
   \nonumber
   &+&\varepsilon\left(\frac{g^{n+1}+ a^T \cdot G}{\Delta t} +  b^T \cdot L(G)\right) = - c^T \cdot {\mum} G- c_{-1} \mu^{n+1} g^{n+1}.
\eea
Multiplying the above equation by the collision invariants $\phi(v)$ and integrating in $v$ gives the moment system
\begin{equation}
\label{eq:IM_NS} \\
\displaystyle   \langle \phi M[f^{n+1}]\rangle = -\sum_{j=0}^{s-1} a_j\,\langle\phi M[f^{n-j}] \rangle - \Delta t\, \sum_{j=0}^{s-1} b_j\left( {\rm div}_x \langle v\otimes \phi M[f^{n-j}]\rangle+\varepsilon \langle \phi\, L(g^{n-j}) \rangle\right).
\end{equation}
Consistency of (\ref{eq:IM_NS}) with the Navier-Stokes equations is obtained if the vector of initial steps satisfies (see \cite{BPESBGK, Golse}) for $j=0,\ldots,s-1$
\bea
\nonumber
\langle \phi\, L(g^{n-j}) \rangle &=& -\left\langle \phi\, L\left(\frac1{\mu}\left(\frac{\partial M[f]}{\partial t}+L(M[f])\right)\right)\Big|_{t=t^{n-j}} \right\rangle+ O(\varepsilon+\Delta t^{q}),\qquad q \geq 1.
\label{eq:IM_WP0}
\eea 
Now we have \cite{Golse}
\be
\left(\frac{\partial M[f]}{\partial t}+L(M[f])\right) = M[f]\left(A(V)\sigma(u)+2B(V)\cdot\nabla_x \sqrt{T}\right)+O(\e),
\label{eq:NSbasic}
\ee
where 
\be
V=({v-{u}})/{\sqrt{{T}}},\quad A(V)=V\otimes V -\frac13|V|^2 I,\quad B(V)=\frac12(|V|^2-5)V,
\label{eq:AB}
\ee
with $\sigma(\cdot)$ defined as in (\ref{eq:NSsigma}).
In (\ref{eq:NSbasic}) we used the dependence on the Maxwellian $M[f]$ on the macroscopic quantities and the Euler system (\ref{eq:Euler}) to replace the time derivatives with space derivatives of the moments.

We assume
\be
g^{n-j}=-\frac1{\mu^{n-j}}M[f^{n-j}]\left(A(V^{n-j})\sigma(u^{n-j})+2B(V^{n-j})\cdot\nabla_x \sqrt{T^{n-j}}\right)+ O(\varepsilon+\Delta t^{q}).
\label{eq:IM_WP}
\ee
Inserting this in (\ref{eq:IM_NS}) gives 
\be
{U^{n+1}=-\sum_{j=0}^{s-1} a_j\, {U}^{n-j}} -  \Delta t\,\sum_{j=0}^{s-1} b_j\left({\rm div}_x{\cal F}(U)^{n-j}
- \varepsilon \, {\rm div}_x {\cal D}(\nabla_x U)^{n-j}\right)+O(\varepsilon^2\Delta t+\varepsilon\Delta t^{q+1}),
\label{eq:NS_limit}
\ee
where 
\[
{\cal D}(\nabla_x {U})^{n-j} = \left(
\begin{array}{c}
 0 \\
 \nu^{n-j}\sigma(u^{n-j}) \\
 \kappa^{n-j}\nabla_x T^{n-j} + \nu^{n-s+1} \sigma(u^{n-j})\cdot u^{n-j} 
\end{array}
\right).  
\]
We omit the detailed computation of the viscosity $\nu^{n-j}$ and heat conduction $\kappa^{n-j}$ in terms of the macroscopic temperature. We refer to \cite{BPESBGK,Golse} for more details. Since the explicit scheme in (\ref{eq:NS_limit}) is an $O(\Delta t^p)$ method for the Navier-Stokes system (\ref{eq:NavierStokes}), one obtains a consistent discretization of (\ref{eq:NavierStokes}) when $\varepsilon\Delta t^p+\Delta t^p=o(\varepsilon)$. 

In order for the scheme to be consistent with the Navier-Stokes limit in time it is crucial that at the next time step $g^{n+1}$ satisfies a relation analogous to (\ref{eq:IM_WP0}). From (\ref{eq:IM_CE}) we can write
\bea
\nonumber
 g^{n+1}&=&  -\frac1{\mu^{n+1}}\left[\frac{c^T}{c_{-1} } \cdot {\mum} G
 +\frac1{c_{-1}}\left(\frac{M[f^{n+1}]+ a^T \cdot M[F]}{\Delta t} +  b^T \cdot L(M[F])\right)\right]\\[-.2cm]
 \label{eq:gn1}
 \\[-.2cm]
 \nonumber
 &&+O\left(\frac{\varepsilon}{\varepsilon+\mu^{n+1}c_{-1}\Delta t}\right).
\eea
For IMEX-BDF schemes, we have $c^T\equiv 0$, and moreover
\be
\frac1{c_{-1}}\left(\frac{M[f^{n+1}]+ a^T \cdot M[F]}{\Delta t} +  b^T \cdot L(M[F])\right) = \left(\frac{\partial M[f]}{\partial t}+L(M[f])\right) \Big|_{t=t^{n+1}} + O(\Delta t^{p}),
\label{eq:IMEX-N1}
\ee
which as a consequence of (\ref{eq:NSbasic}) yields the desired  result
\bea
\nonumber
 g^{n+1}&=& -\frac1{\mu^{n+1}}M[f^{n+1}]\left(A(V^{n+1})\sigma(u^{n+1})+2B(V^{n+1})\cdot\nabla_x \sqrt{T^{n+1}}\right) \\[-.2cm]
 \label{eq:gn1b}
 \\[-.2cm]
 \nonumber
 &&+ O\left(\varepsilon+\Delta t^p+\frac{\varepsilon}{\varepsilon+\mu^{n+1}c_{-1}\Delta t}\right).
\eea

For more general IMEX linear multistep methods high order accuracy when estimating the time derivative is lost, and an estimate similar to (\ref{eq:IMEX-N1}) but with first order accuracy in time holds true. Therefore, thanks to the order conditions (\ref{eq:GIMEXcond}), from (\ref{eq:NSbasic}) and (\ref{eq:gn1}) we get
\bea
\nonumber
 g^{n+1}&=& -\frac1{\mu^{n+1}}M[f^{n+1}]\left(A(V^{n+1})\sigma(u^{n+1})+2B(V^{n+1})\cdot\nabla_x \sqrt{T^{n+1}}\right) \\[-.2cm]
\label{eq:gn1c}
 \\[-.2cm]
 \nonumber
 &&+ O\left(\varepsilon+\Delta t+\frac{\varepsilon}{\varepsilon+\mu^{n+1}c_{-1}\Delta t}\right).
\eea


We have proved the following:
\begin{theorem}
If the vector of initial steps is well-prepared with respect to the Navier-Stokes limit as defined in (\ref{eq:IM_WP}), then, for small values of $\varepsilon$ and with $\varepsilon\Delta t^q+\Delta t^p=o(\varepsilon)$, the IMEX multistep scheme (\ref{eq:IMEX_BGK}) becomes the explicit multistep scheme 
(\ref{eq:NS_limit}) for the Navier-Stokes system (\ref{eq:NavierStokes}), with $q=p$ for the IMEX-BDF methods and $q=1$ for the other IMEX multistep methods.
\end{theorem}

\begin{remark}
\begin{itemize}
\item The result shows that the schemes are capable, in principle, to capture the Navier-Stokes asymptotics without resolving the small scale $\varepsilon$.
Note, however, that when $\Delta t$ is of order $\varepsilon$, the error term in (\ref{eq:gn1b}) and (\ref{eq:gn1c}) is an $O(1/(1+\mu^{n+1}c_{-1}))$ term. Therefore, we may expect loss of accuracy of the schemes in such a regime.
\item
 It is well-known that the limiting Navier-Stokes system for the BGK model differs from the classical one derived from the full Boltzmann equation since it corresponds to a Prandtl number $Pr=1$, whereas the classical one for a monoatomic gas is $Pr=2/3$. The correct Prandtl number can be recovered considering a velocity dependent collision frequency $\mu=\mu(x,v,t)$ as in the ES-BGK model \cite{BPESBGK}.
 \item The analysis just performed can be carried on in a similar way also for the full Boltzmann equation. The conclusions one obtains for the various IMEX multistep schemes are exactly the same as for the BGK model, therefore the schemes are capable to describe correctly the Navier-Stokes regime. Here we omit the details. We will, however, discuss the details of the Navier-Stokes asymptotics for the full Boltzmann equation using a penalized approach in the next section.     
 \end{itemize}
\end{remark}

\section{Penalized IMEX multistep schemes for the Boltzmann
equation} 
In the case of the full Boltzmann equation, although the approach just described remains formally valid, a major difficulty concerns the need to solve the system of
nonlinear equations originated by the application of an implicit method
to the collision operator $Q_B(f,f)$ introduced in (\ref{eq:Q}). The computational cost of such
integral operator, characterized by a five fold nonlinear integral
which depends on the seven dimensional space $(x,v,t)$, is extremely
high and makes it extremely expensive the use of iterative solvers. It is however interesting to remark that for IMEX multistep methods since a single new evaluation of the collision operator is required at each time step, a single inversion would be also required which in principle is more feasible then the case of IMEX Runge-Kutta schemes.

To overcome this difficulty, the idea is to reformulate the collision
part using a suitable penalization term. This idea, using a BGK
model as penalization term, has been introduced recently
in~\cite{Filbet}. We refer to~\cite{dimarco6,dimarco7,Qin} for  extensions of
this approach to high-order IMEX Runge-Kutta methods and exponential Runge-Kutta methods.

\subsection{Asymptotic preserving penalized IMEX linear multistep schemes} 
Let us now denote with $Q_{P}(f)$ a general operator which will be
used to penalize the original Boltzmann operator $Q_{B}(f,f)$. The characteristics of $Q_P(f)$ are to be computable and
invertible at a low computational cost and that it preserves the
local equilibrium, namely $Q_P(f)=0$ implies $f=M[f]$.

We will then rewrite the collision operator in the form
\be
Q_B(f,f)=(Q_B(f,f)-Q_P(f))+Q_P(f)=G_P(f,f)+Q_P(f),
\label{eq:QP}
\ee
where by construction $\langle\phi\, G_P(f,f)\rangle=0$, and the corresponding kinetic equation reads
\be
\partial_t
f+v\cdot\nabla_x
f=\frac{1}{\varepsilon}G_P(f,f)+\frac{1}{\varepsilon}Q_P(f).\label{eq:Bref2}
\ee
Now, we consider IMEX multistep schemes in which only
the simpler operator $Q_P(f)$ is treated implicitly, while the
term $G_P(f,f)$ describing the deviations of the true Boltzmann operator $Q_B(f,f)$ from the simplified operator $Q_P(f)$ and the
convection term $\nabla_x f$ are treated explicitly. 

We can now
introduce the general class of penalized IMEX linear multistep schemes for the
Boltzmann equation in the form

\bea
\nonumber
   f^{n+1} &=& - \sum_{j=0}^{\nuu-1} a_j f^{n-j} + \Delta t \sum_{j=0}^{\nuu-1} b_j\, \left(\frac1{\varepsilon}G_P(f^{n-j},f^{n-j})-v\cdot\nabla_x f^{n-j}\right)\\[-.25cm]
   \label{eq:GIMEXb}
   \\[-.25cm]
   \nonumber
   &&+ \Delta t \sum_{j=-1}^{\nuu-1} c_j \frac{1}{\varepsilon}Q_P(f^{n-j}),
\eea
or equivalently, using vector notation as
\begin{equation}
   f^{n+1} = - a^T \cdot F + \Delta t\, b^T \cdot \left(\frac1{\varepsilon}G_P(F,F) -L(F)\right) + \frac{\Delta t}{\varepsilon} c^T \cdot Q_P(F)+\frac{\Delta t}{\varepsilon} c_{-1} Q_P(f^{n+1}),
\label{eq:GIMEXbv}
\end{equation}
where $G_P(F,F)=(G_P(f^{n},f^n),\ldots,G_P(f^{n-\nuu+1},f^{n-\nuu+1}))^T$. 

\begin{remark}
The above penalized approach is equivalent to start from the standard IMEX multistep scheme (\ref{eq:GIMEXv}) and to introduce a penalization only on the implicit term. More precisely we can write (\ref{eq:GIMEXv}) as
\begin{equation}
   f^{n+1} = - a^T \cdot F - \Delta t\, b^T \cdot L(F) + \frac{\Delta t}{\varepsilon} c^T \cdot Q_B(F,F)+\frac{\Delta t}{\varepsilon} c_{-1} G_P(f^{n+1},f^{n+1})+\frac{\Delta t}{\varepsilon} c_{-1} Q_P(f^{n+1}).
\label{eq:GIMEXv3}
\end{equation}
Now thanks to the order conditions (\ref{eq:GIMEXcond}) we have 
\be
\frac{b^T-c^T}{c_1{-1}}G_P(F,F)=G_P(f^{n+1},f^{n+1})+O(\Delta t^p).
\label{eq:lagrange}
\ee
In fact, by direct inspection, one observes that 
\be
\frac{b_i-c_i}{c_{-1}}=\prod_{k=0 \atop k\neq j}^{s-1}\frac{k+1}{k-j},
\ee
are the Lagrange interpolation weights at time $n+1$. Using (\ref{eq:lagrange}) in (\ref{eq:GIMEXv3}) we recover (\ref{eq:GIMEXbv}).
\label{rk:Lagrange}
\end{remark}

We want now to derive conditions for the penalized IMEX multistep schemes to be AP and
asymptotically accurate. We have the following result.
\begin{theorem}
\label{th:pap} If the vector of initial steps is well-prepared, in the limit
$\varepsilon\rightarrow 0$, scheme
(\ref{eq:GIMEXbv}) becomes the explicit multistep scheme characterized by
($a, b$) applied to the limit
Euler system (\ref{eq:Euler}).
\end{theorem}
 
\proof The proof is similar to that of Theorem \ref{th:ap}. Let us first note that if we multiply the IMEX
multistep scheme (\ref{eq:GIMEXbv}) by the collision
invariants $\phi(v)=1, v, v^2$ and integrate the result in velocity
space we obtain the same explicit multistep method (\ref {eq:GIMEXfm}) for the
moment system (\ref{eq:macr}).

Now we can write equation (\ref{eq:GIMEXbv})
in the form 
\begin{equation}
   \varepsilon f^{n+1} = - \varepsilon a^T \cdot F + \Delta t\, b^T \cdot \left(G_P(F,F) -\varepsilon L(F)\right)+\Delta t\, c^T \cdot Q_P(F) +\Delta t\, c_{-1} Q_P(f^{n+1}),
   \label{eq:espiIMEXb}
\end{equation}
which as $\varepsilon\to 0$ yields
\begin{equation}
   0 =  b^T \cdot G_P(F,F) + c^T \cdot Q_P(F) + c_{-1} Q_P(f^{n+1}).
   \label{eq:espiIMEXb2}
\end{equation}
Since the initial steps are well-prepared, as $\varepsilon\rightarrow 0$ we get $f^{n-j}=M[f^{n-j}]$, $j=0,\ldots,\nuu-1$ which implies $Q_P(F)\equiv 0$, $G_P(F,F)\equiv 0$ and therefore since $c_{-1}\neq 0$ we have
\be
Q_P(f^{n+1})=0 \Rightarrow f^{n+1}=M[f^{n+1}].
\ee
Thus (\ref{eq:GIMEXfm})
 becomes the explicit multistep method (\ref{eq:GIMEXfeuler}) for to the limiting Euler system (\ref{eq:Euler}).\\
$\Box$

Note that, for IMEX-BDF schemes in the case of arbitrary initial steps, equation (\ref{eq:espiIMEXb2}) reduces to
\begin{equation}
   0 =  b^T \cdot G_P(F,F) + c_{-1} Q_P(f^{n+1}).
   \label{eq:espiIMEXbdf}
\end{equation}
Therefore, at variance with the non penalized case, the scheme does not project $f^{n+1}$ over its local equilibrium $M[f^{n+1}]$ unless the initial vector is well-prepared so that $G_P(F,F)\equiv 0$.

\subsection{Penalized IMEX relaxation multistep schemes for the Boltzmann equation}
The choice of the optimal penalization operator depends essentially on a balance between accuracy and computational cost. A popular choice which optimize efficiency is given by the simple BGK relaxation operator \cite{dimarco7, Filbet}. In fact, with this choice, a fundamental property of equations
(\ref{eq:GIMEXbv}) with $Q_P(f)=\mu(M[f]-f)$ is that they can be solved explicitly.

The penalized IMEX multistep method takes the form
\bea
\nonumber
   f^{n+1} &=& - \sum_{j=0}^{\nuu-1} a_j f^{n-j} + \Delta t \sum_{j=0}^{\nuu-1} b_j\, \left(\frac1{\varepsilon}G_P(f^{n-j},f^{n-j})-v\cdot\nabla_x f^{n-j}\right)\\[-.25cm]
   \label{eq:PIMEXBGK}
  \\[-.25cm]
  \nonumber
 &&  + \Delta t \sum_{j=-1}^{\nuu-1} c_j \frac{\mu^{n-j}}{\varepsilon}(M[f^{n-j}]-f^{n-j}),
\eea
or equivalently, using vector notation as
\bea
\nonumber
   f^{n+1} &=& - a^T \cdot F + \Delta t\, b^T \cdot \left(\frac1{\varepsilon}G_P(F,F) -L(F)\right) \\[-.25cm]
   \\[-.25cm]
   \nonumber
   &&+ \frac{\Delta t}{\varepsilon} \left(c^T \cdot {\mum}(M[F]-F)+ c_{-1}\mu^{n+1}(M[f^{n+1}]-f^{n+1})\right),
\eea
where the only implicit term is the factor $\mu^{n+1}(M[f^{n+1}]-f^{n+1})$ in which $M[f^{n+1}]$ and $\mu^{n+1}$ depend only on the moments $\langle\phi f^{n+1}\rangle$. If we now integrate the above equation agains the collision invariants thanks to the conservations (\ref{eq:QC}) we obtain again the moment scheme (\ref{eq:GIMEXfm}). Thus $\langle \phi f^{n+1}\rangle$, and so $M[f^{n+1}]$ and $\mu^{n+1}$, can be explicitly evaluated and system (\ref{eq:GIMEXfm}) is explicitly solvable. 

\subsubsection{The space homogeneous case}
Next we focus on the monotonicity properties of the penalized IMEX schemes. As a prototype problem we restrict to the space homogeneous case for the simplified BGK model, which in this case however involves the full penalized IMEX method. As we will see, the presence of the explicit scheme due to the introduction of the penalization lead to different (but still very severe) non negativity conditions.

The penalized IMEX method for $Q_B(f,f)=\eta(M[f]-f)$ in the homogeneous case reads
\bea
\nonumber
   f^{n+1} &=& - a^T \cdot F + \frac{\Delta t}{\e}\, b^T \cdot \left((\eta-\mu)(M[f^n]\,e-F)\right) \\[-.25cm]
   \\[-.25cm]
   \nonumber
   &&+ \frac{\Delta t}{\varepsilon} \left(c^T \cdot {\mu}(M[f]\,e-F)+ c_{-1}\mu(M[f^n]-f^{n+1})\right),
\eea
 or equivalently
\be
   f^{n+1} 
   = -\left(\lambda a^T+(1-\lambda) \left(\frac{c^T+\xi b^T}{c_{-1}}\right)\right)\cdot F
   +(1-\lambda)\left(\frac{c^T+\xi b^T}{c_{-1}} \cdot e + 1\right) M[f^n].
      \label{eq:s1p}
   \ee   
where $\xi=(\eta-\mu)/\mu \in (-1,\infty)$ is the penalization factor.   

Setting $z=\mu\Delta t/\varepsilon$ we have $\lambda=1/(1+c_{-1}z)$ and we can state 
\begin{proposition}
Sufficient conditions to guarantee that $f^{n+1}\geq 0$ when $F\geq 0$ in (\ref{eq:s1p}) are that
\begin{eqnarray}
\label{eq:posnp}
a^T +z \left(c^T + \xi b^T\right) \leq 0,\\
(c^T\cdot e + c_{-1})(1+\xi) \geq 0.
\label{eq:posn1p}
\end{eqnarray}
\end{proposition}
Again conditions (\ref{eq:posnp})-(\ref{eq:posn1p}) must be
interpreted component by component and for $\xi=0$ reduce to (\ref{eq:posn})-(\ref{eq:posn1}).  In particular, (\ref{eq:posnp})
depend on $z$ and originates a time step restriction (and constraints over the choice of the penalization factor $\xi$). Note that in (\ref{eq:posn1p}) we used the fact that $b^T \cdot e=c^T\cdot e+c_{-1}$. For penalized IMEX-BDF schemes $c^T = 0$ and conditions (\ref{eq:posnp}) give 
\[
\max_{0\leq j\leq s-1} \{a_j/(\xi |b_j|), b_j < 0\} \leq z \leq \min_{0\leq j\leq s-1} \{-a_j/(\xi b_j), b_j > 0\},\qquad \xi > 0. 
\]
Note that, for schemes of order higher then $2$ the above condition can never be satisfied.
In the case of penalized IMEX-Adams methods conditions (\ref{eq:posnp}) become 
\[
z\leq 1/(c_0+\xi b_0), \quad c_0+\xi b_0>0,\qquad 
\max_{1\leq j\leq s-1}\{c_j/|b_j|, b_j < 0\} \leq \xi \leq \min_{1\leq j\leq s-1}\{-c_j/b_j, b_j > 0\}.
\]

For example, using notations of Table \ref{tb:examples}, nonnegativity time step restrictions of some penalized IMEX multistep schemes are reported in Table \ref{tb:timestep}.

\begin{table}[t]
\caption{Nonnegativity restrictions of some penalized multistep IMEX schemes}
\begin{center}
{\small
\begin{tabular}{l|c c}
\hline\\[-.25cm]
Scheme & Time step & Penalization\\
\hline\\[-.25cm]
IMEX-BDF1 & $z \leq 1/\xi$ & $\xi>0$\\[+.25cm]
IMEX-CN2 & $z \leq 2/(1+3\xi)$ & $\xi \geq 0$\\[+.25cm]
IMEX-MCN2 & $z \leq 8/(3+12\xi)$ & $\xi\geq 1/8$\\[+.25cm]
IMEX-BDF2 & $1/(2\xi)\leq z \leq 1/\xi$ & $\xi>0$\\[+.25cm]
IMEX-SG2 & $z \leq \min\{1/2,1/(2\xi)\}$ & $\xi\geq 0$\\[+.25cm]
IMEX-AD3 & $z \leq 12/(1551/2500+23\xi)$ & $107/2196\leq\xi\leq 492/4147$\\[+.25cm]
\hline
\end{tabular}
}
\end{center}
\label{tb:timestep}
\end{table}


Let us remark that since from the first order conditions $a^T\cdot e = -1$ we can also write
\be
f^{n+1}=-\frac{a^T+z(c^T+\xi b^T)}{1+c_{-1}z}\cdot F+\left(\frac{e^T}{\nuu}+\frac{a^T+z(c^T+\xi b^T)}{1+c_{-1}z}\right)\cdot eM.\label{eq:simpp2}\ee
Therefore we can state
\begin{proposition}
If conditions (\ref{eq:posnp})-(\ref{eq:posn1p}) are satisfied 
then (\ref{eq:simpp2}) is a convex combination of the initial steps and the Maxwellian state.
\end{proposition}

From the above proposition we obtain the following entropy bound
\be
H(f^{n+1}) \leq \sum_{j=0}^{\nuu-1} \alpha(\xi)_j H(f^{n-j})+\left(1-\sum_{j=0}^{\nuu-1}\alpha(\xi)_j\right) H(M) \leq \max_{0\leq j \leq s-1} H(f^{n-j}),
\ee
where now
\[
\alpha(\xi)_j = -\frac{a_j+z(c_j+\xi b_j)}{1+c_{-1}z}.
\]


\begin{figure}
\begin{center}
\includegraphics[scale=0.35]{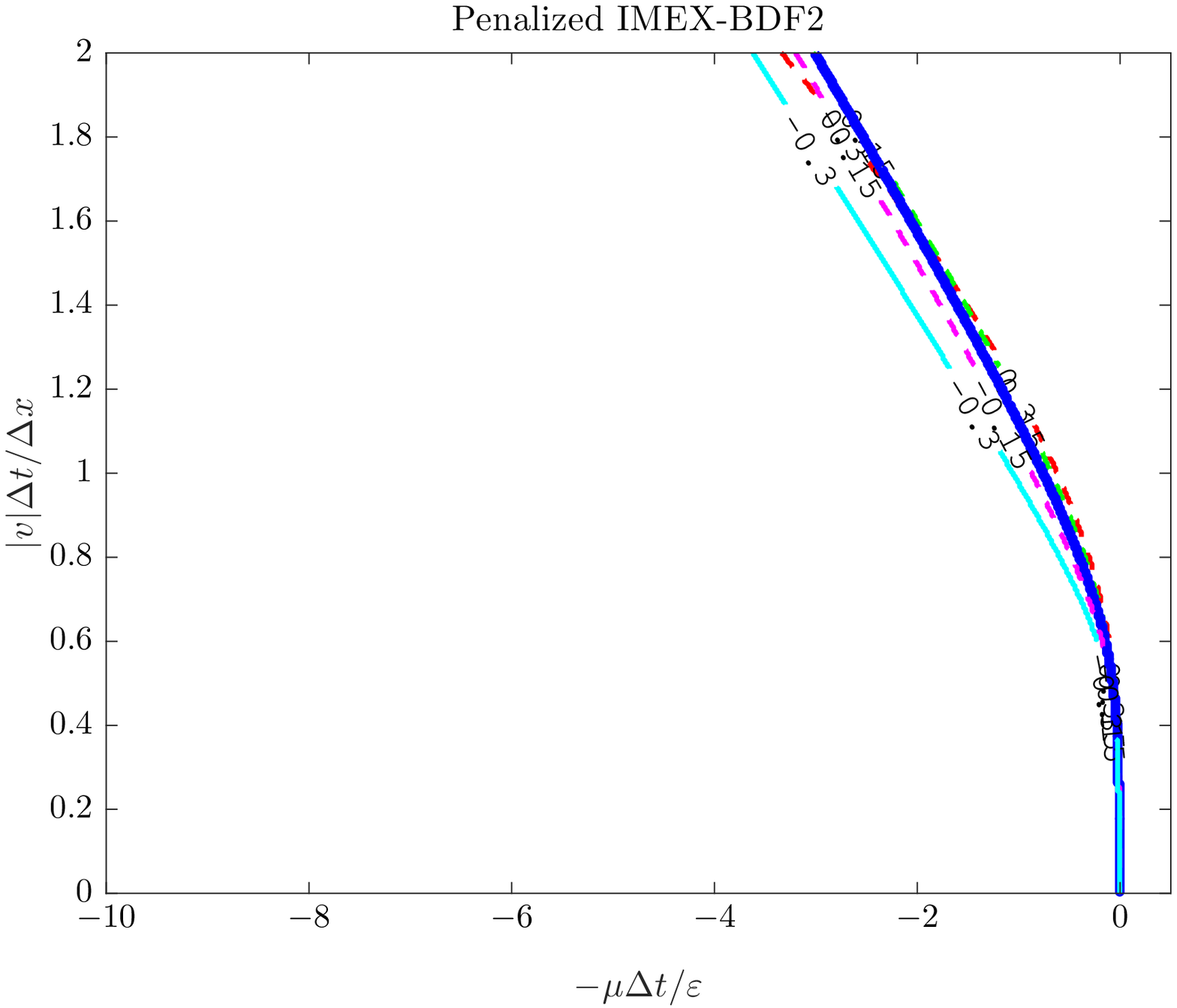}\qquad
\includegraphics[scale=0.35]{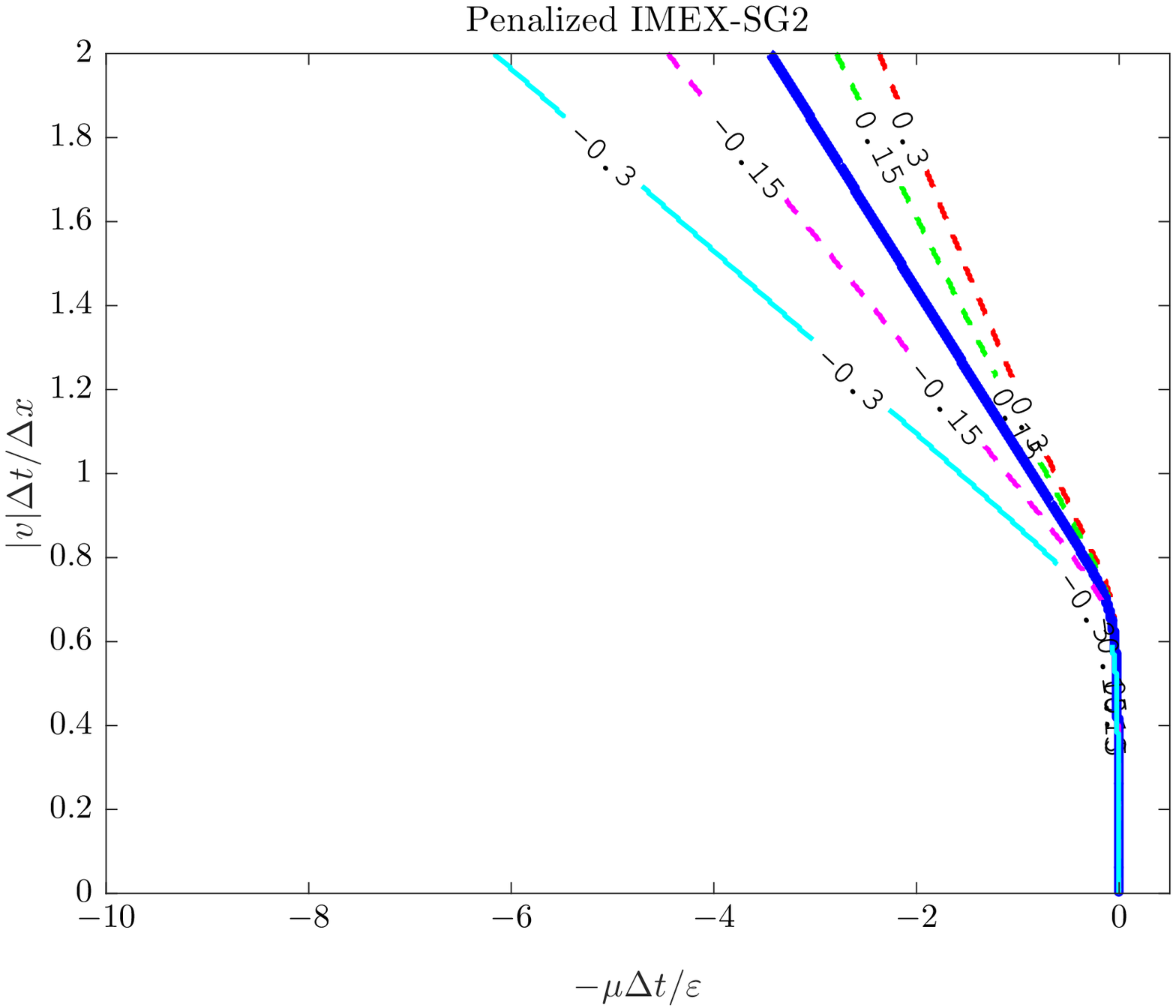}\\[+.1cm]
\includegraphics[scale=0.35]{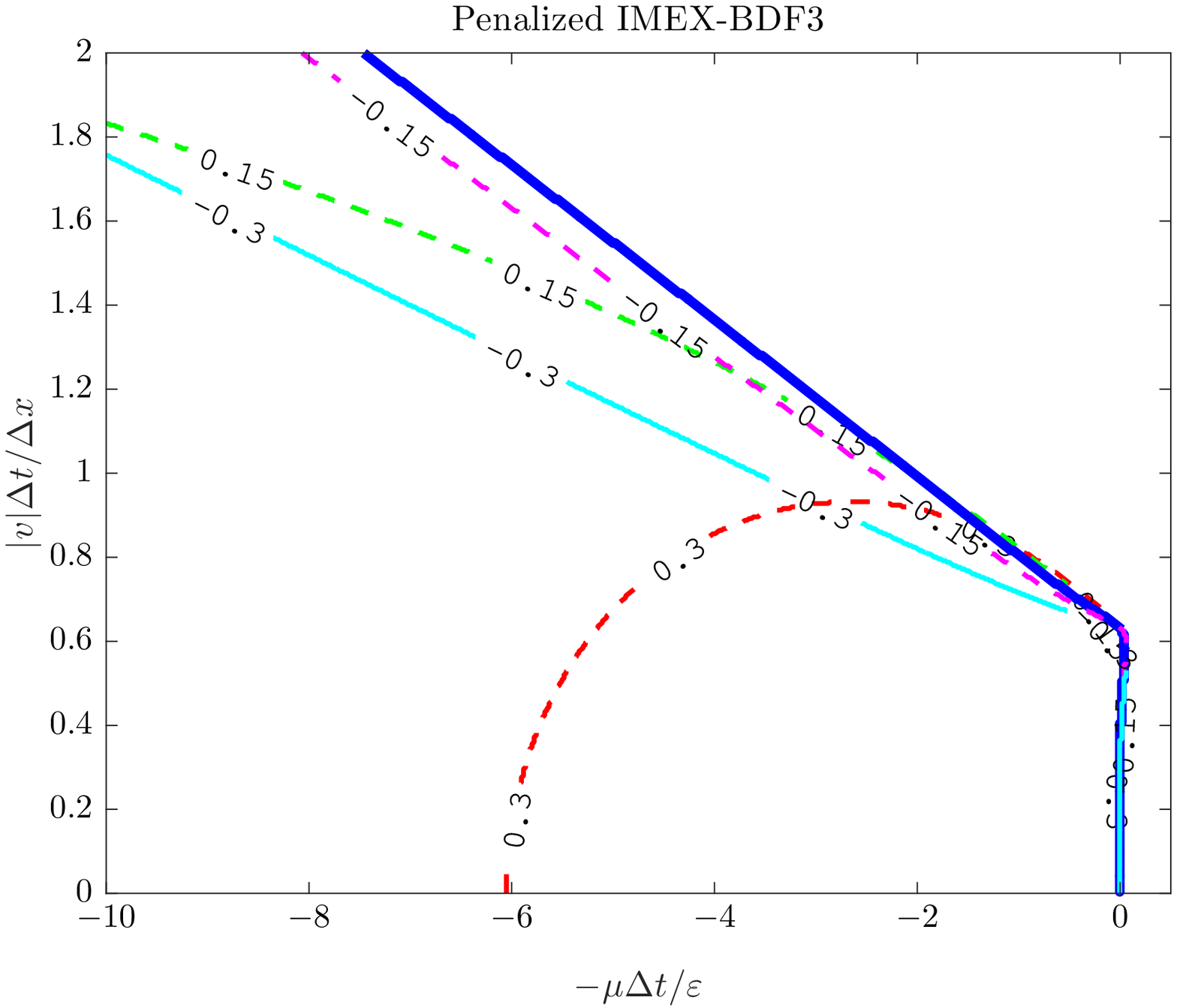}\qquad
\includegraphics[scale=0.35]{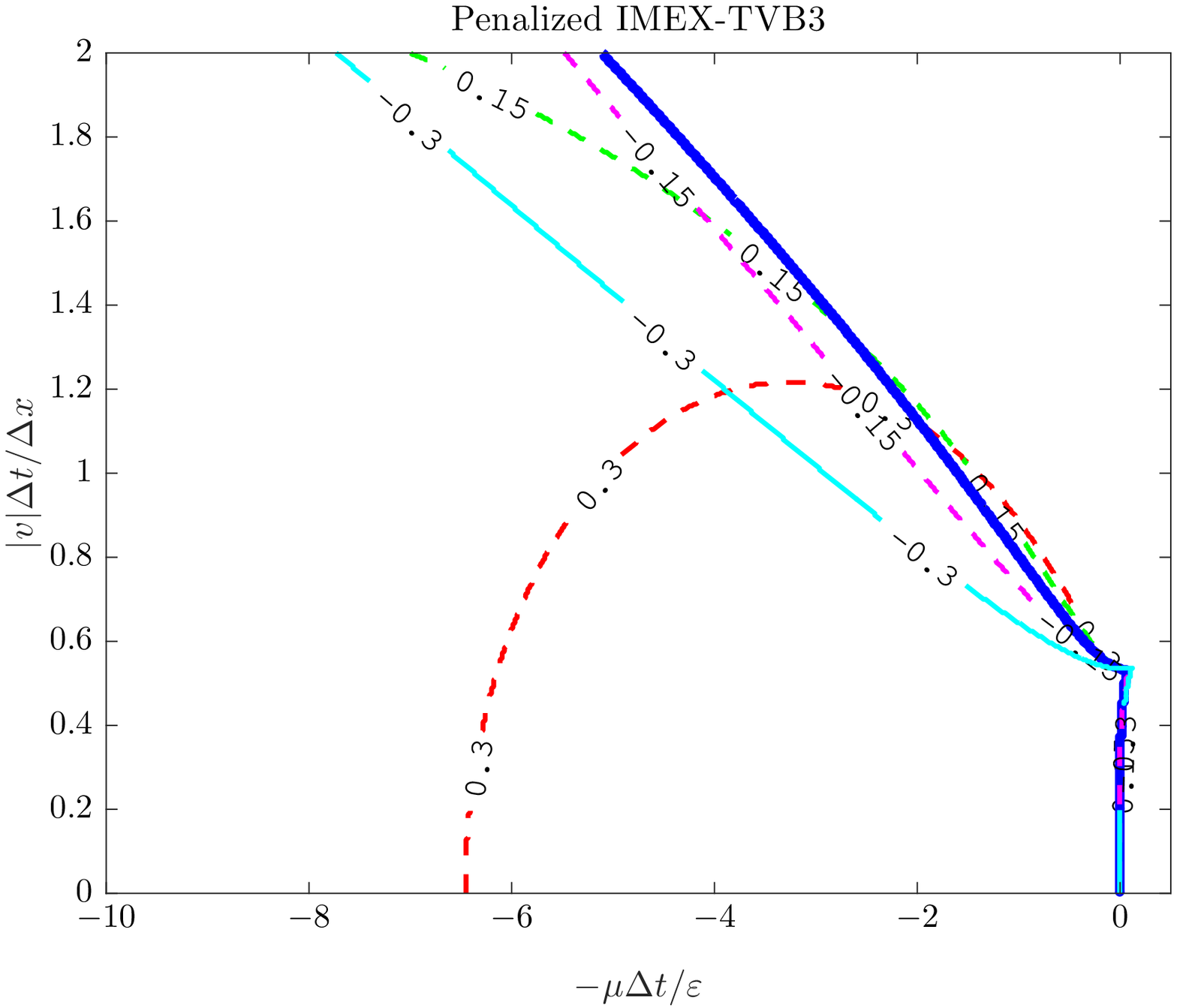}
\end{center}
\caption{Boundary of the stability region for problem (\ref{eq:ebgk})-(\ref{eq:cd}). Penalized second order (top) and third order (bottom) IMEX multistep methods with a penalization factor $\xi=-0.3,-0.15,0,0.15,0.3$. The label on the contour lines denotes the penalization factor. The thick contour line without label correspond to $\xi=0$.}
\label{fig:stab1}
\end{figure}

\begin{figure}
\begin{center}
\includegraphics[scale=0.35]{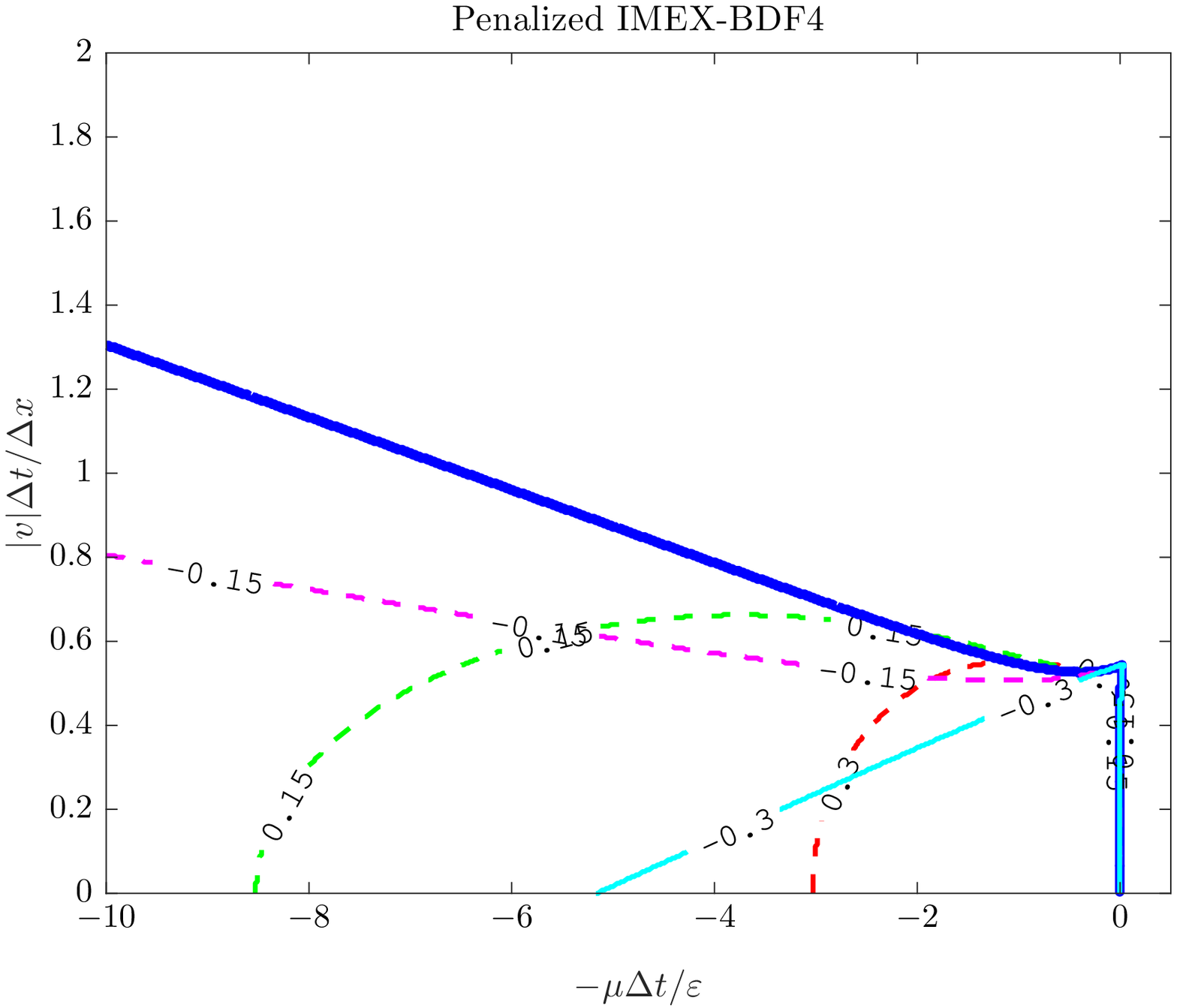}\qquad
\includegraphics[scale=0.35]{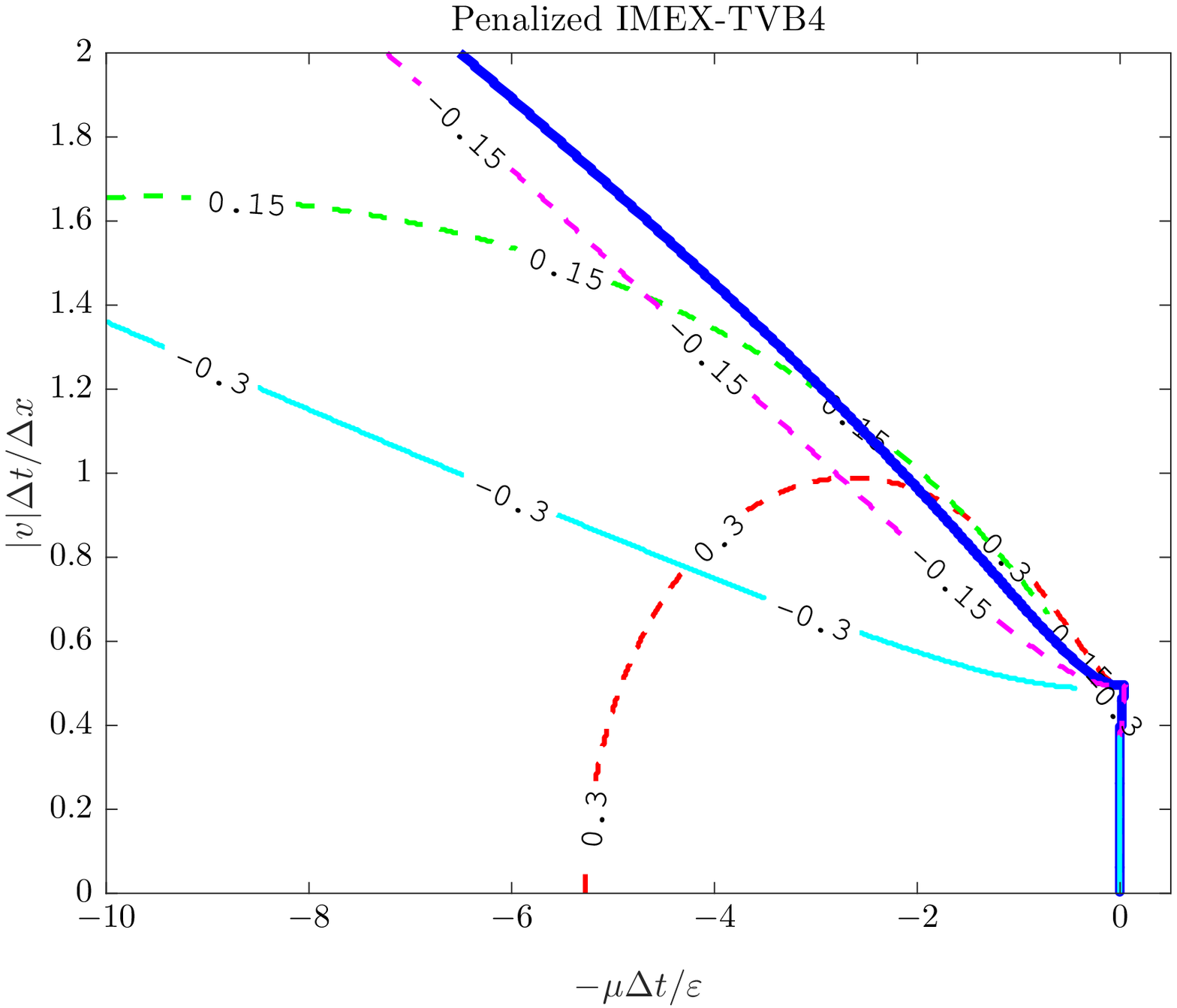}\\[+.1cm]
\includegraphics[scale=0.35]{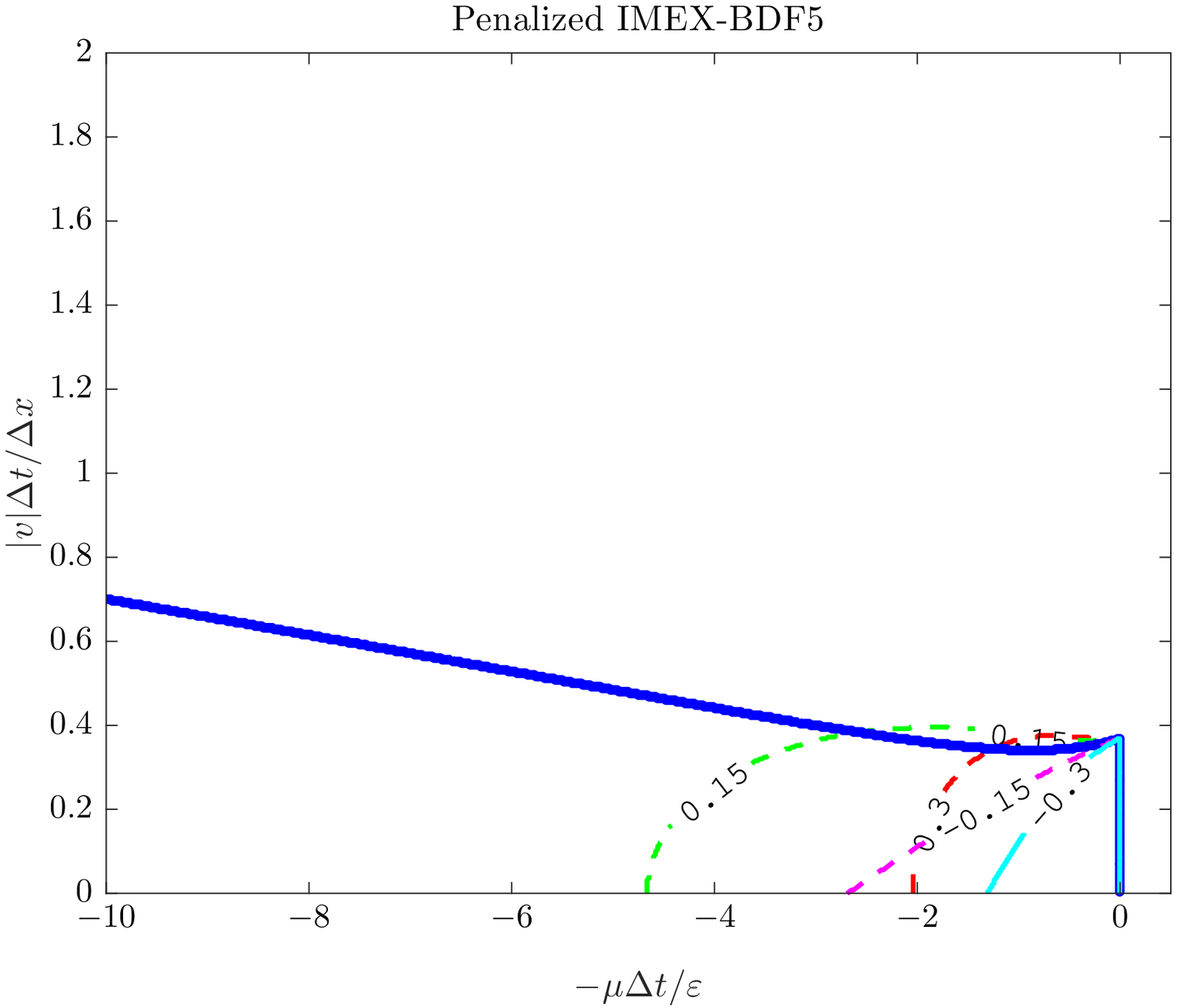}\qquad
\includegraphics[scale=0.35]{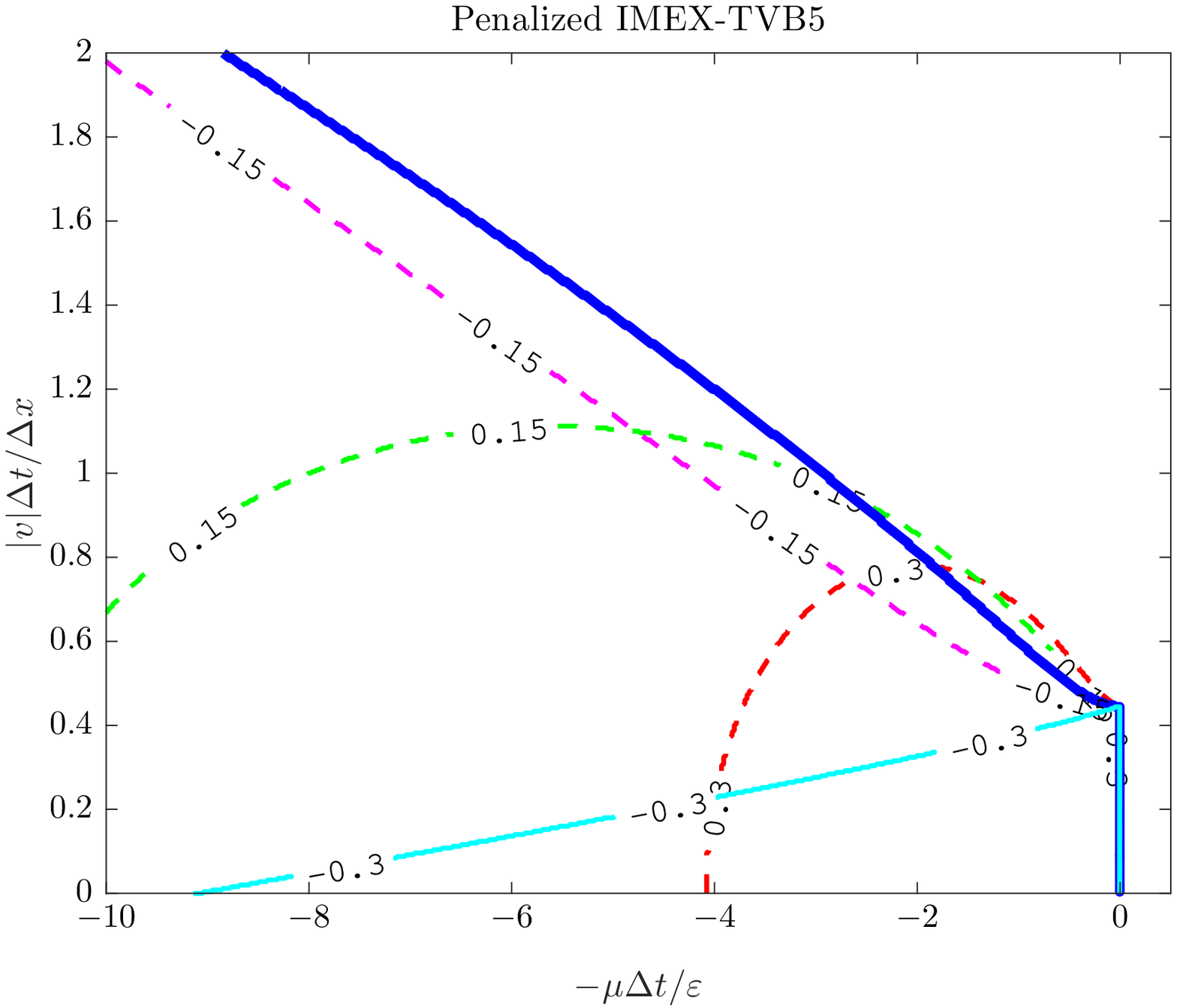}
\end{center}
\caption{Boundary of the stability region for problem (\ref{eq:ebgk})-(\ref{eq:cd}). Penalized fourth order (top) and fifth order (bottom) IMEX multistep methods with a penalization factor $\xi=-0.3,-0.15,0,0.15,0.3$. The label on the contour lines denotes the penalization factor. The thick contour line without label correspond to $\xi=0$.}
\label{fig:stab2}
\end{figure}

\subsubsection{Stability restrictions}
\label{sec:stab}
It is clear that the penalized approach may reduce the stability region of the full IMEX scheme applied without penalization. This can be observed using the test equation
\be
y'=\lambda_E y + \lambda_I y,\quad \lambda_E,\lambda_I \in \C
\label{eq:test}
\ee
where the term $\lambda_E y$ is evaluated explicitly and the term $\lambda_I y$ implicitly. Typically, for prototype kinetic equations, we can assume $\lambda_E$ pure imaginary and $\lambda_I$ real and non positive. In fact, if we consider a simple transport equation with source like
\be
f_t + v f_x = -\frac{\eta}{\varepsilon}f,
\label{eq:ebgk}
\ee
after discretization of the space derivatives, by Fourier transform we recover the test equation (\ref{eq:test}) with $\lambda_E$ the eigenvalues of the advection term and $\lambda_I$ those of the stiff source term. We have $\lambda_I=-\eta/\varepsilon$ and, for example using central differences for (\ref{eq:ebgk}), we get  
\be
\lambda_E=iv\sin(2k)/\Delta x, 
\label{eq:cd}
\ee
where $k$ is the frequency of the corresponding Fourier mode. 

Applying the IMEX multistep method gives the characteristic equation
\be
\rho(\zeta)-z_E\sigma_E(\zeta)-z_I\sigma_I(\zeta)=0,
\label{eq:char}
\ee
with $z_E=\lambda_E \Delta t$ and $z_I=\lambda_I \Delta t$ and 
\[
\rho(\zeta)=\zeta^s+a^T\cdot{\cal X},\quad \sigma_E(\zeta)=b^T\cdot{\cal X},\quad  \sigma_I(\zeta)=c_{-1}\zeta^s+c^T\cdot{\cal X}, 
\]
where ${{\cal X}}=(\zeta^{s-1},\zeta^{s-2},\ldots,\zeta,1)^T$.
Stability corresponds to the requirement that all roots of (\ref{eq:char}) have modulus less or equal one and that all multiple roots have modulus less then one. Even if the prototype equation is very simple and its analysis, in general, cannot be extended to linear systems (obtained for example trough the methods of lines), for practical purposes it is often used to obtain informations for stability \cite{Ascher, Ascher2, HR, HRS}. 

For the penalized IMEX approach 
\be
y'=\lambda_E y + (\lambda_I-\lambda_P) y + \lambda_P y,
\ee
with $\lambda_P$ real and non positive, we obtain the same structure of the characteristic equation (\ref{eq:char}) with the only difference that now $z_E=(\lambda_E+\lambda_I-\lambda_P) \Delta t$ and $z_I=\lambda_P \Delta t$. In the case of (\ref{eq:ebgk}) solved by central differences, if $\lambda_P=-\mu/\varepsilon$, by introducing the penalization factor $\xi=(\eta-\mu)/\mu$ we have $z_E=(iv\sin(2k)/\Delta x-\xi\mu/\varepsilon)\Delta t$ and $z_I=-\mu \Delta t/\varepsilon$.
The stability regions of various penalized IMEX multistep schemes are reported in Figures \ref{fig:stab1} and \ref{fig:stab2}. The boundaries separate the stability region (lower part) from the instability region (upper part).
Analogous results for some non penalized IMEX multistep methods and some IMEX Runge-Kutta methods have been presented in \cite{Ascher, Ascher2, HR}. Its is clear that high order methods are more sensitive to the choice of the penalization factor and may lead to instabilities unless the eigenvalues of the stiff part are estimated with enough accuracy. This is particularly true for fourth and fifth order IMEX-BDF methods.

\subsubsection{Navier-Stokes asymptotics}

Let us rewrite the penalized IMEX scheme as
\bea
   \nonumber
   \frac{f^{n+1}+ a^T \cdot F}{\Delta t} +  b^T \cdot L(F) &=&  \frac1{\e} b^T \cdot Q_B(F,F)\\[-.25cm]
   \label{eq:PIMEXBGK2}
   \\[-.2cm]
   \nonumber
   &&+\frac{1}{\varepsilon} (c^T-b^T) \cdot {\mum}(M[F]-F)+\frac{1}{\varepsilon} c_{-1} \mu^{n+1}(M[f^{n+1}]-f^{n+1}),
\eea
and consider the discrete Chapman-Enskog expansion taking
\be
f^{n+1}=M[f^{n+1}]+\varepsilon g^{n+1},\quad F=M[F]+\varepsilon G,
\ee
where, using the same notation of section 3, we have $\langle \phi g^{n+1}\rangle=0$ and $\langle \phi G\rangle \equiv 0$. 

Inserting the above expansions in the numerical method yields
\bea
\nonumber
   \frac{M[f^{n+1}]+ a^T \cdot M[F]}{\Delta t} &+&  b^T \cdot L(M[F])+\varepsilon\left(\frac{g^{n+1}+ a^T \cdot G}{\Delta t} +  b^T \cdot L(G)\right) \\[-.25cm]
   \label{eq:IM_CEb}
   \\[-.25cm]
   \nonumber
   &=& b^T\cdot{\cal L}_M(G)+\e b^T Q(G)- (c^T-b^T) \cdot {\mum}G-  c_{-1} \mu^{n+1} g^{n+1},
\eea
where ${\cal L}_M(G)=({\cal L}_{M[f^n]}(g^n),\ldots,{\cal L}_{M[f^{n-s+1}]}(g^{n-s+1}))^T$ and 
\be
{\cal L}_{M[f^{n-j}]}(g^{n-j})=2 Q_B(M[f^{n-j}],g^{n-j}),\quad j=0,\ldots,s-1
\label{eq:LM}
\ee
is the linearized Boltzmann operator with respect to $M[f^{n-j}]$.

Multiplying the above equation by the collision invariants $\phi(v)$ and integrating in $v$ gives the moment system
\begin{equation}
\label{eq:IM_CEb2} \\
\displaystyle   \langle \phi M[f^{n+1}]\rangle = -\sum_{j=0}^{s-1} a_j\,\langle\phi M[f^{n-j}] \rangle - \Delta t\, \sum_{j=0}^{s-1} b_j\left( {\rm div}_x \langle v\otimes \phi M[f^{n-j}]\rangle+\varepsilon \langle \phi\, v\cdot\nabla_x g^{n-j} \rangle\right).
\end{equation}
Consistency with the Navier-Stokes equations corresponds to assume that the vector of initial values $G$ satisfies for $j=0,\ldots,s-1$
\bea
\nonumber
\langle \phi\, v\cdot\nabla_x g^{n-j} \rangle &=& -\left\langle \phi\, v\cdot\nabla_x\left({\cal L}_{M[f]}^{-1}\left(\frac{\partial M[f]}{\partial t}+{v}\cdot \nabla_x M[f]\right)\right){\Big |_{t=t^{n-j}}}\right\rangle+ O(\varepsilon)+O(\Delta t^q), \quad q \geq 1.
\label{eq:IM_WP1}
\eea 
This is guaranteed if \cite{Caflisch,cercignani,Golse} 
\be
g^{n-j}=-{\cal L}_{M[f]}^{-1}\left(\frac{\partial M[f]}{\partial t}+{v}\cdot \nabla_x M[f]\right)\Big|_{t=t^{n-j}} + O(\varepsilon+\Delta t^q).
\label{eq:NSconst}
\ee
Using (\ref{eq:NSbasic}) and the properties of the linearized collision operator we get
\bea
\nonumber
&&{\cal L}_{M[f]}^{-1}\left(\frac{\partial M[f]}{\partial t}+{v}\cdot \nabla_x M[f]\right) + O(\varepsilon)\\
\nonumber
&=&{\cal L}_{M[f]}^{-1}\left(M[f]\left(A(V)\sigma(u)+2B(V)\cdot\nabla_x \sqrt{T}\right)\right)+ O(\varepsilon),\\
\label{eq:IM_WP2}
&=& -{\cal L}_{M[f]}^{-1}\left(M[f]A(V)\right)\sigma(u)-2{\cal L}_{M[f]}^{-1}\left(M[f]B(V)\right)\cdot\nabla_x \sqrt{T}+ O(\varepsilon),\\
\nonumber
&=& -\frac{M[f]}{\rho}\left(a(T,|V|)A(V)\sigma(u)+2b(T,|V|) B(V)\cdot\nabla_x \sqrt{T}\right)+ O(\varepsilon),
\eea
where $V=({v-{u}})/{\sqrt{{T}}}$, $A(V)$, $B(V)$ have been defined in (\ref{eq:AB}), and $a(T,|V|)$, $b(T,|V|)$ are suitable scalar functions.

Therefore, if we assume
\bea
\nonumber
g^{n-j}&=&-\frac{M[f^{n-j}]}{\rho^{n-j}}\left(a(T^{n-j},|V^{n-j}|)A(V^{n-j})\sigma(u^{n-j})\right.\\[-.25cm]
\label{eq:NSassum}
\\[-.25cm]
\nonumber
&&\left.+2b(T^{n-j},|V^{n-j}|) B(V^{n-j})\cdot\nabla_x \sqrt{T^{n-j}}\right)+ O(\varepsilon+\Delta t^q),
\eea
we have that (\ref{eq:IM_CEb2}) corresponds to
\be
{U^{n+1}=-\sum_{j=0}^{s-1} a_j\, {U}^{n-j}} -  \Delta t\,\sum_{j=0}^{s-1} b_j\left({\rm div}_x{\cal F}(U)^{n-j}
- \varepsilon \, {\rm div}_x {\cal D}(\nabla_x U)^{n-j}\right)+O(\varepsilon^2\Delta t)+O(\varepsilon\Delta t^{q+1}),
\label{eq:NScorrect}
\ee
where 
\[
{\cal D}(\nabla_x {U})^{n-j} = \left(
\begin{array}{c}
 0 \\
 \nu^{n-j}\sigma(u^{n-j}) \\
 \kappa^{n-j}\nabla_x T^{n-j} + \nu^{n-s+1} \sigma(u^{n-j})\cdot u^{n-j} 
\end{array}
\right).  
\]
We omit the detailed computation of the viscosity $\nu^{n-j}$ and heat conduction $\kappa^{n-j}$ in terms of $a(T^{n-j},|V^{n-j})$ and $b(T^{n-j},|V^{n-j})$. We refer to \cite{Golse} for more details. Since the discretization in (\ref{eq:NScorrect}) is an $O(\Delta t^p)$ approximation of the Navier-Stoke system (\ref{eq:NavierStokes}) we obtains the correct Navier-Stokes limit if $\varepsilon\Delta t^q+\Delta t^p = o(\varepsilon)$.

However, to achieve order $p$ consistency at the next time step, $g^{n+1}$ should satisfy the analogous of (\ref{eq:NSconst}). 
From (\ref{eq:IM_CEb}) we have
\bea
\nonumber
g^{n+1} = -\frac1{\mu^{n+1}}&&\left[\frac{(c^T-b^T)}{c_{-1}}\cdot{\mum}G+\frac1{c_{-1}}\left(
   \frac{M[f^{n+1}]+ a^T \cdot M[F]}{\Delta t} +  b^T \cdot L(M[F])\right)\right.\\[-.2cm]
   \\[-.2cm]
   \nonumber
    &&-\frac{b^T}{c_{-1}}\cdot{\cal L}_M(G)\Big]+O\left(\frac{\varepsilon}{\varepsilon+c_{-1}\Delta t}\right).
\eea
Now using (\ref{eq:NSassum}) and Remark \ref{rk:Lagrange}, we obtain
\[
\frac{(c^T-b^T)}{c_{-1}}\cdot{\mum}G = \mu^{n+1}{\cal L}^{-1}_{M[f]}\left(\frac{\partial M[f]}{\partial t}+{v}\cdot \nabla_x M[f]\right)\Big|_{t=t^{n+1}}+O(\varepsilon+\Delta t^p+\Delta t^q).
\]
Moreover, for IMEX-BDF schemes we have
\[
\frac1{c_{-1}}\left(
   \frac{M[f^{n+1}]+ a^T \cdot M[F]}{\Delta t} +  b^T \cdot L(M[F])\right)=\left(\frac{\partial M[f]}{\partial t}+{v}\cdot \nabla_x M[f]\right)\Big|_{t=t^{n+1}}+O(\Delta t^p)
\]
and since $c^T=0$
\[
\frac{b^T}{c_{-1}}\cdot{\cal L}_M(G)=\left(\frac{\partial M[f]}{\partial t}+{v}\cdot \nabla_x M[f]\right)\Big|_{t=t^{n+1}}+O(\Delta t^p).
\]
Finally we get 
\be
g^{n+1}=-{\cal L}_{M[f]}^{-1}\left(\frac{\partial M[f]}{\partial t}+{v}\cdot \nabla_x M[f]\right)\Big|_{t=t^{n+1}} + O(\varepsilon+\Delta t^q+\Delta t^p)+O\left(\frac{\varepsilon}{\varepsilon+c_{-1}\Delta t}\right).
\label{eq:NSconst2}
\ee
For more general IMEX multistep methods, high order accuracy in the evaluation of the time derivatives is lost and we obtain the same estimate (\ref{eq:NSconst2}) but with first order accuracy in time. This shows that the penalized IMEX multistep schemes are able to capture correctly the Navier-Stokes asymptotics. 

We can state the following:
\begin{theorem}
If the vector of initial steps is well-prepared with respect to the Navier-Stokes limit as defined in (\ref{eq:IM_WP1}), then, for small values of $\varepsilon$ and with $\varepsilon\Delta t^q+\Delta t^p=o(\varepsilon)$, the penalized IMEX relaxation multistep scheme (\ref{eq:PIMEXBGK2}) becomes the explicit multistep scheme 
(\ref{eq:NScorrect}) for the Navier-Stokes system (\ref{eq:NavierStokes}), with $q=p$ for the IMEX-BDF methods and $q=1$ for the other IMEX multistep methods.
\end{theorem}

\begin{remark}
\begin{itemize}
\item Again, even if the schemes are capable to capture the Navier-Stokes asymptotics without resolving the small scale $\varepsilon$ we may observe loss of accuracy of the schemes when $\Delta t$ is of order $\varepsilon$,since  the error term becomes an $O(1/(1+\mu^{n+1}c_{-1}))$ term. 
 \item The same kind of analysis holds true also for more general penalization operators $Q_P(f)$ and the conclusions one obtains for the various penalized IMEX multistep schemes are the same as in the case of the BGK penalization.     
 \end{itemize}
\end{remark}

\section{Numerical examples}
Scope of this section is to perform several numerical tests for the IMEX multistep schemes
in order to highlight their behaviors in various regimes of the stiffness parameter $\varepsilon$. We will consider both the BGK and the Boltzmann case.
In the case of the full Boltzmann equation the penalized IMEX multistep schemes have been applied using the BGK relaxation operator as penalization term. 
We start our numerical results with the space non homogeneous BGK equation, next we pass to the space homogeneous Boltzmann equation
and we conclude with the full non homogeneous Boltzmann equation. In all tests we have considered eight different IMEX multistep schemes from order two up to order five. More precisely, we report results for the IMEX-BDF schemes
of order two, three, four and five, for the second order IMEX-SG scheme, and for the IMEX-TVB schemes of order three, four and five. The coefficients of the schemes are given in Table \ref{tb:examples}.

\subsection{Non homogeneous BGK equation}
We numerically measure the accuracy of the IMEX multistep methods in solving 
equation (\ref{eq:1b}) with the BGK collision operator (\ref{BGK}) on a periodic smooth
solution. The computation is performed on $(x, v) \in [0, 1] \times
[-v_{\max}, v_{\max}]$, with $v_{\max} = 8$. We take $N_v = 512$
grid points in the velocity space. A 5th order WENO
scheme for the space discretization \cite{Shu} has been used. The largest time step
is fixed equal to $\Delta t_{max} = \frac{\Delta x}{4v_{\max}}$ for all simulations, which is
sufficient to guarantee stability. The initial data is \be
\varrho_{0}(x)=\frac{2+\sin(8\pi x)}{3}, \ u_{0}(x)=0, \ T_{0}(x)=\frac{2+\cos(8 \pi x)}{3},\ee 
where an initial distribution, $f(x,v,t=0)=f_0=M[f_0]+\varepsilon g$ with $g$ a perturbation from equilibrium consistent with the Navier-Stokes limit has been used as in (\ref{eq:IM_WP0}).
We report the convergence rates with respect to the time step $\Delta t$ measured by computing the $L_1$
norm of the error for the density for different values of the
Knudsen number, i.e. $\varepsilon=10^{-1}$, $\varepsilon=10^{-2}$
and $\varepsilon=10^{-5}$. 
In order to perform such a measure of the convergence rate we repeat the computation for a decreasing time step
$\Delta t_{1} = \Delta t_{max}/2$, $\Delta t_{2} = \Delta t_{max}/4$  and $\Delta t_{3} = \Delta t_{max}/8$
while the grid points in space is fixed to $N_x=128$. Finally, in order to initialize the different multistep schemes a third order IMEX Runge-Kutta scheme \cite{dimarco7}
has been employed with very small time steps to produce the vector of initial data accurate up to machine precision.

\begin{figure}
\begin{center}
\includegraphics[scale=0.35]{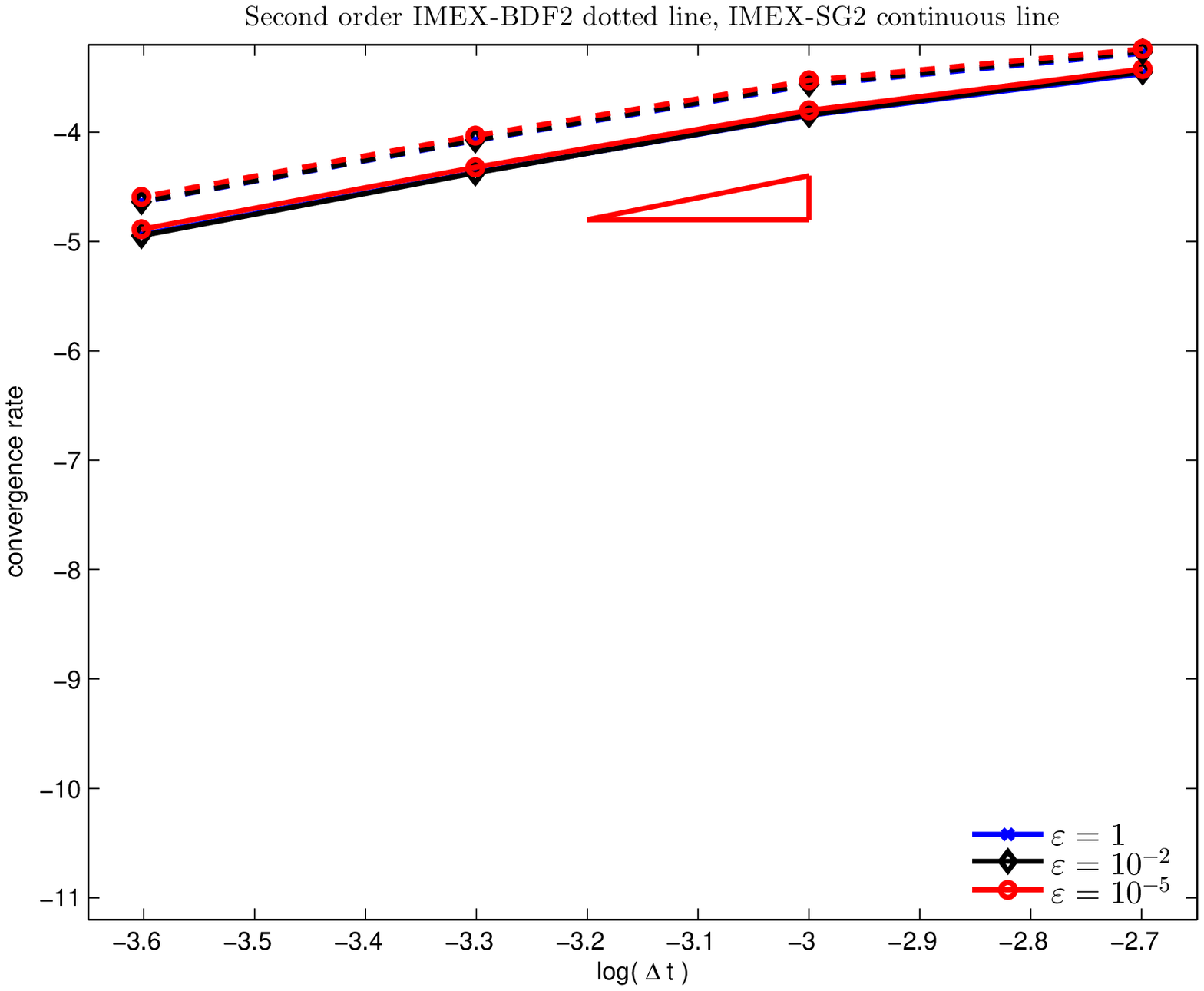}
\includegraphics[scale=0.35]{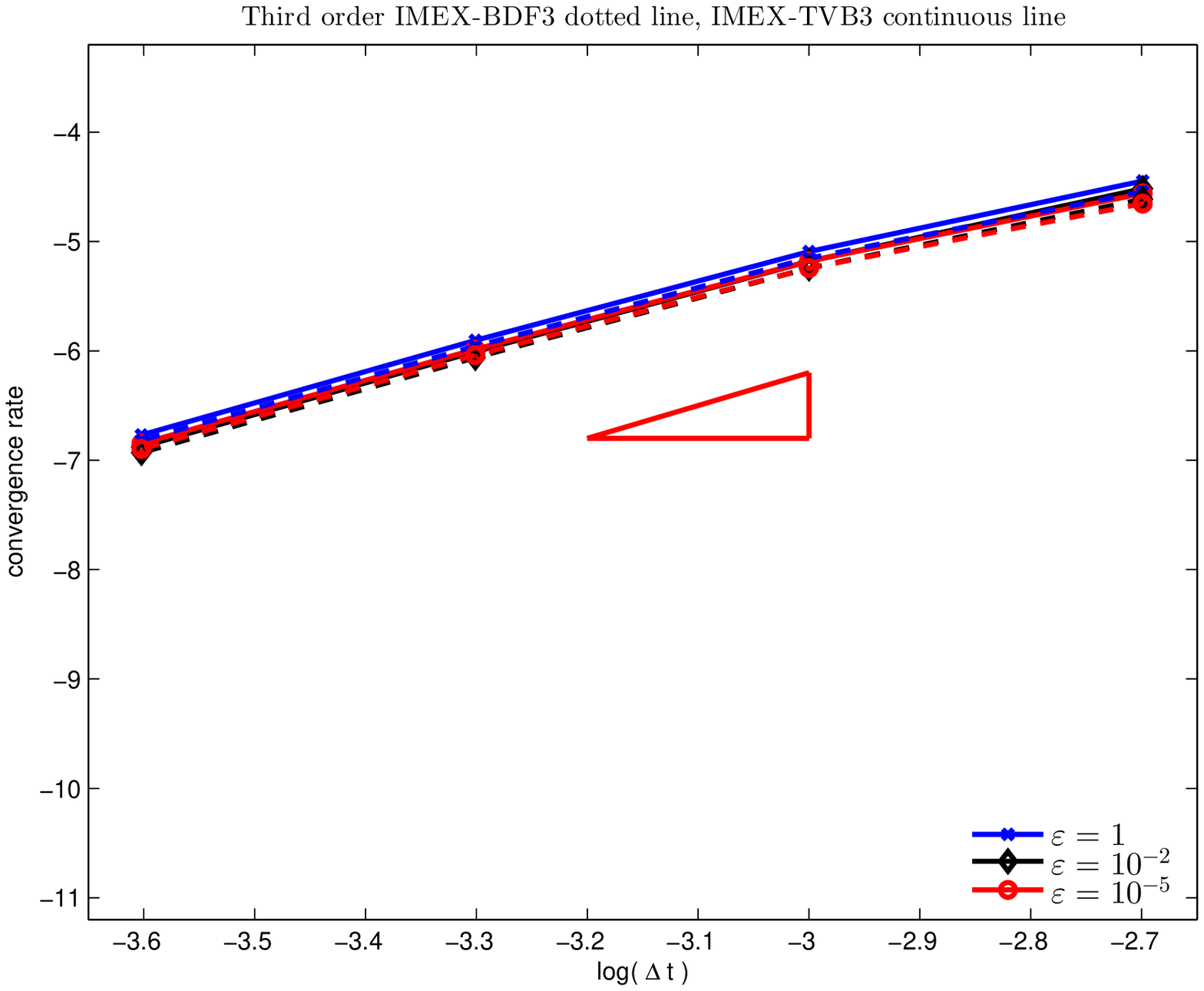}\\
\includegraphics[scale=0.35]{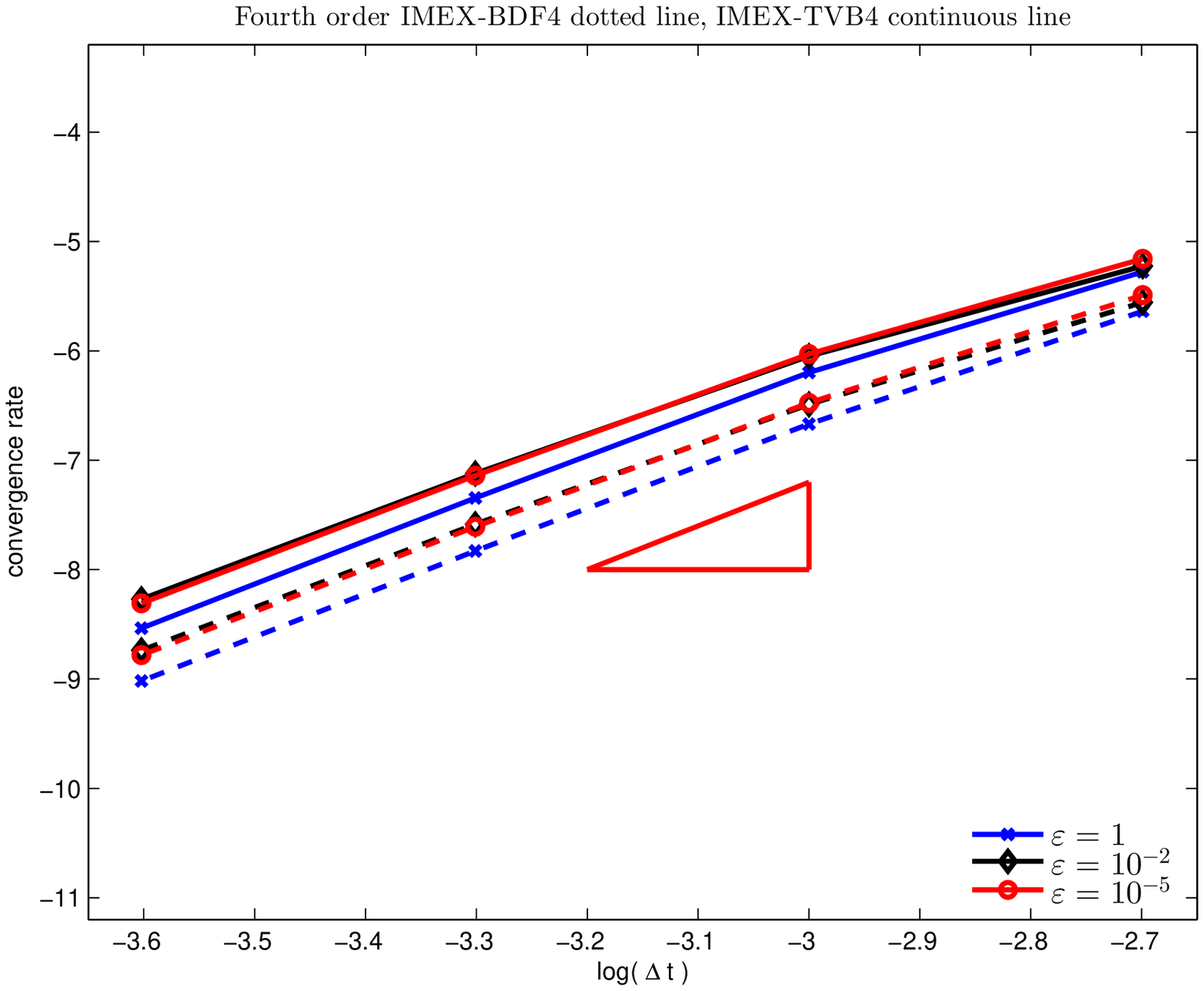}
\includegraphics[scale=0.35]{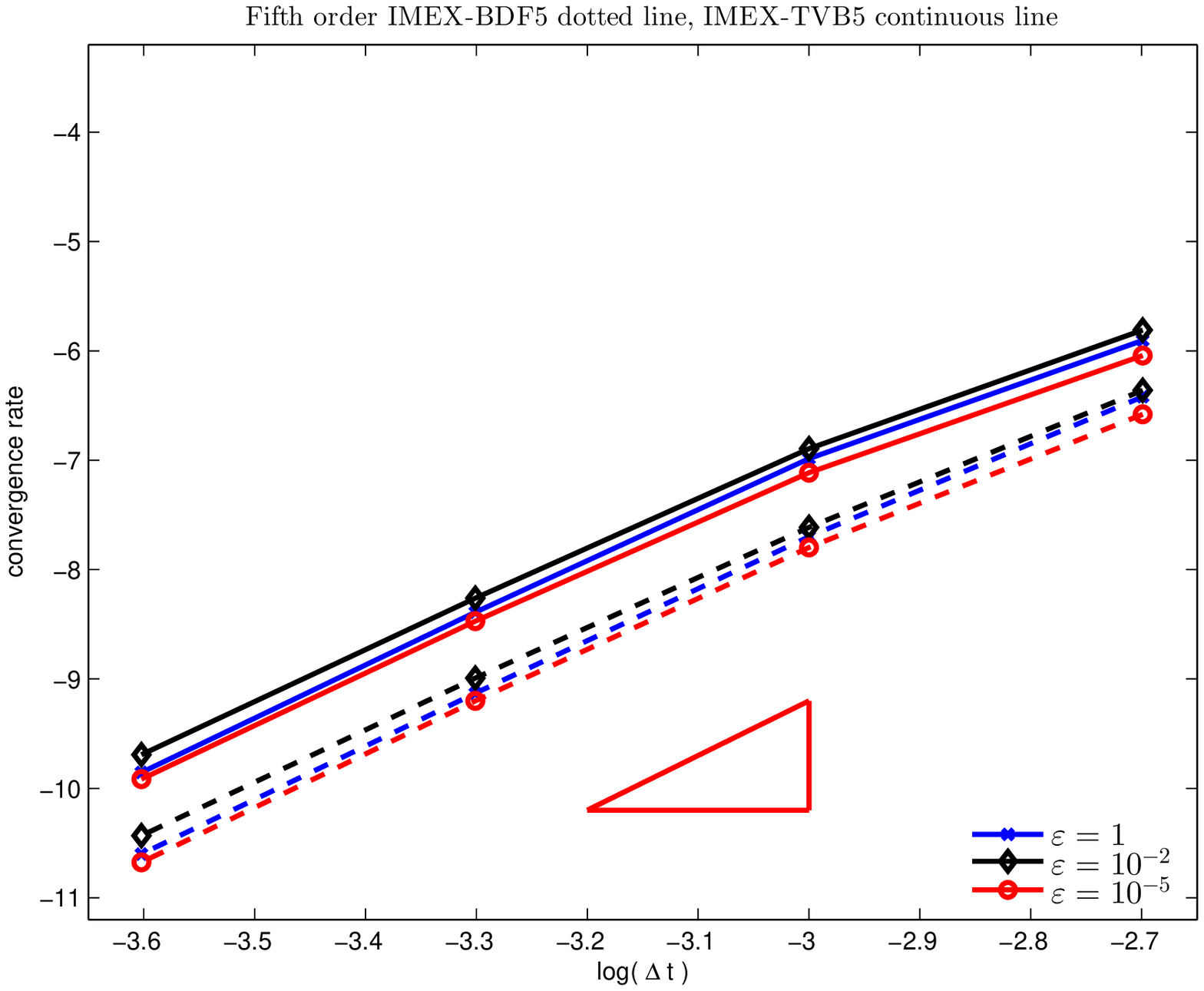}
\caption{$L_{1}$ error for the density $\rho$. Top left IMEX-BDF2 and IMEX-SG2 schemes. Top right IMEX-BDF3 and IMEX-TVB3.
Bottom left IMEX-BDF4 and IMEX-TVB4, bottom right IMEX-BDF5 and IMEX-TVB5.}\label{fig:conv_BGK}
\end{center}
\end{figure}

Figure \ref{fig:conv_BGK} shows the results for all the multistep schemes tested. 
The second order schemes are depicted on the top left of the figure, the third order schemes on the top right,
the fourth order schemes are reported on the bottom left while the fifth order schemes are shown on the bottom right side.
All schemes tested exhibit the theoretical rate of convergence for the different values of $\varepsilon$. We point out that in our numerical simulation we have observed that fourth and fifth order schemes are more sensitive to the correct initialization of the initial vector in agreement with the Navier-Stokes system, otherwise loss of accuracy has been observed. 

 \begin{figure}
 \begin{center}
 \includegraphics[scale=0.35]{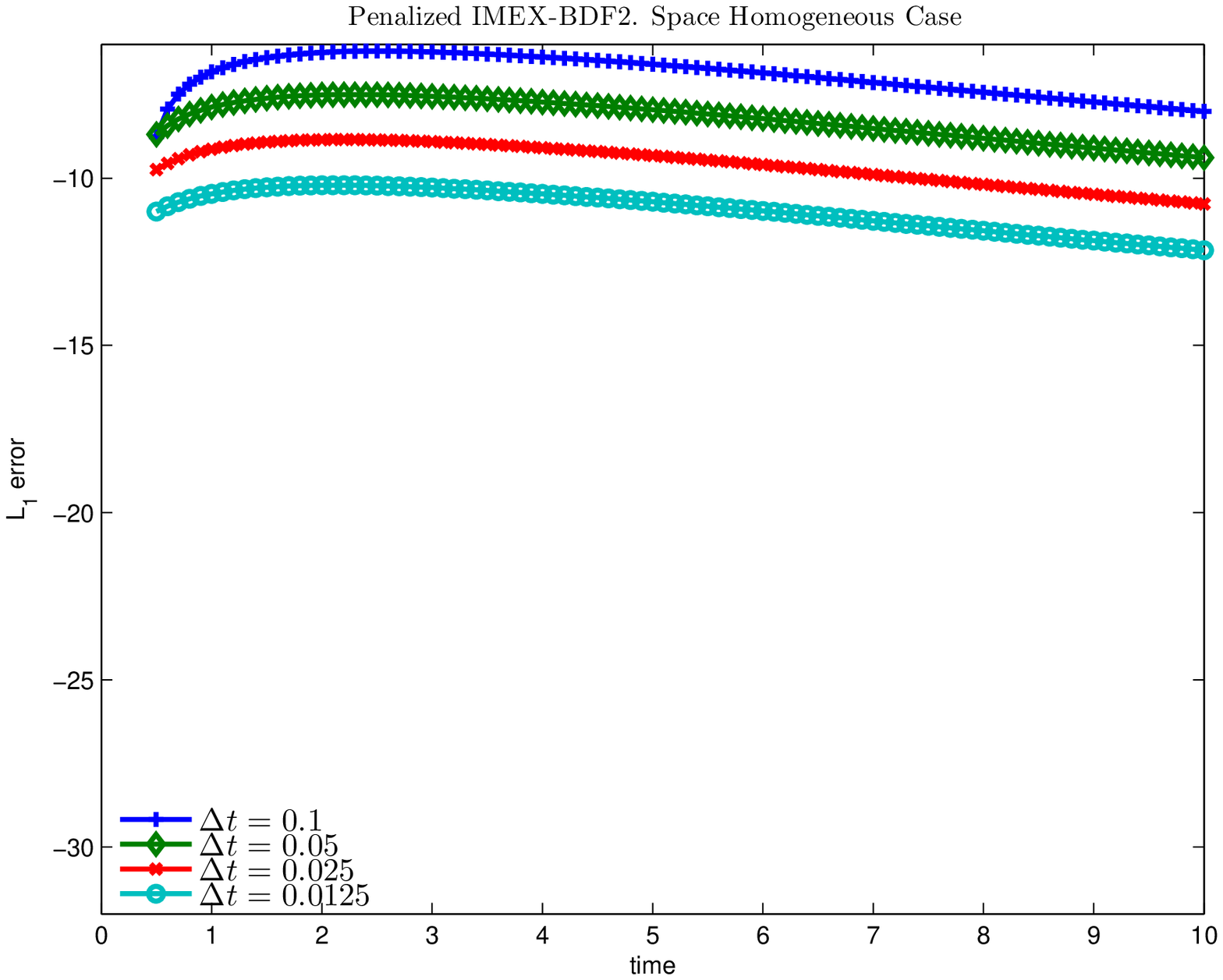}
 \includegraphics[scale=0.35]{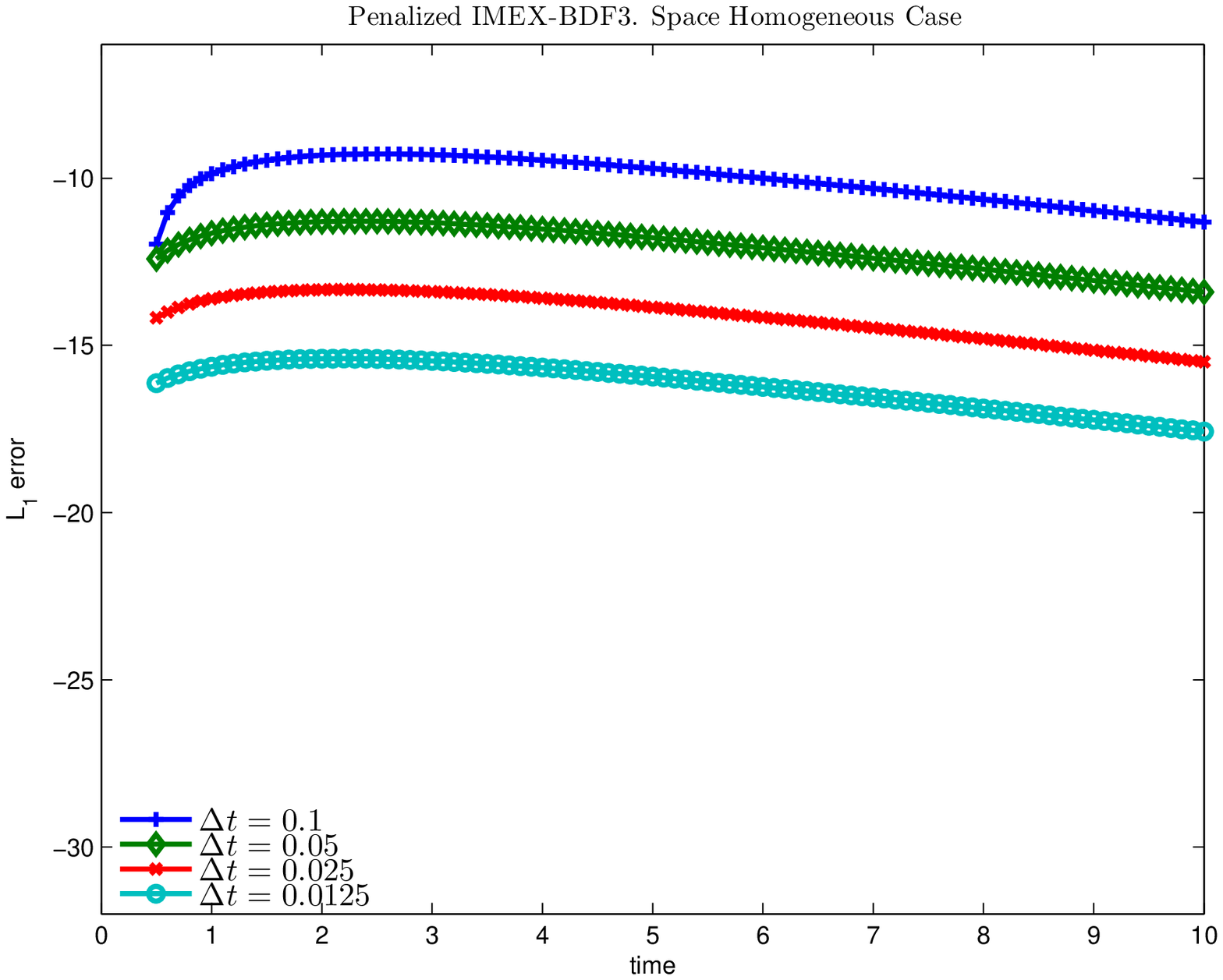}\\
\includegraphics[scale=0.35]{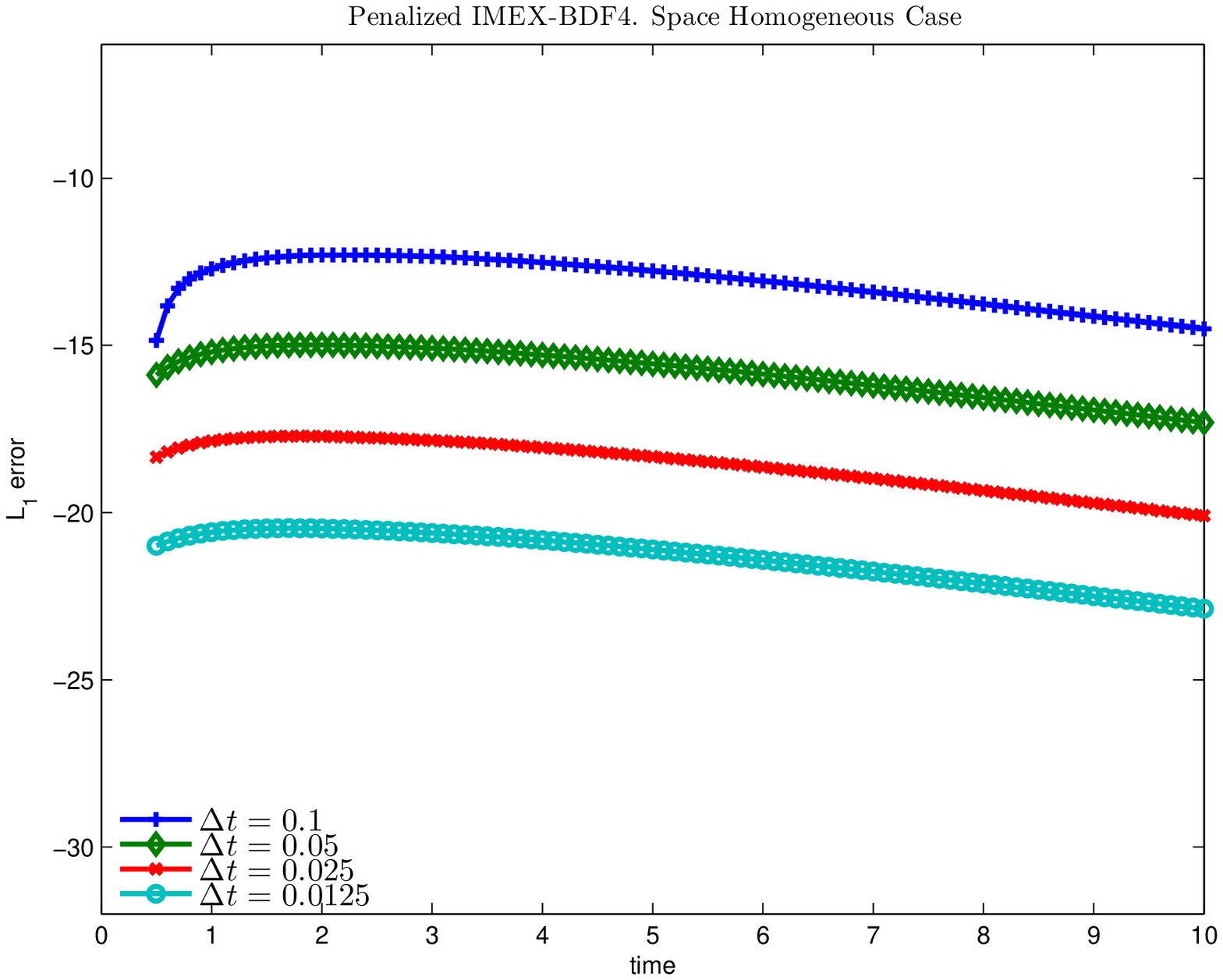}
\includegraphics[scale=0.35]{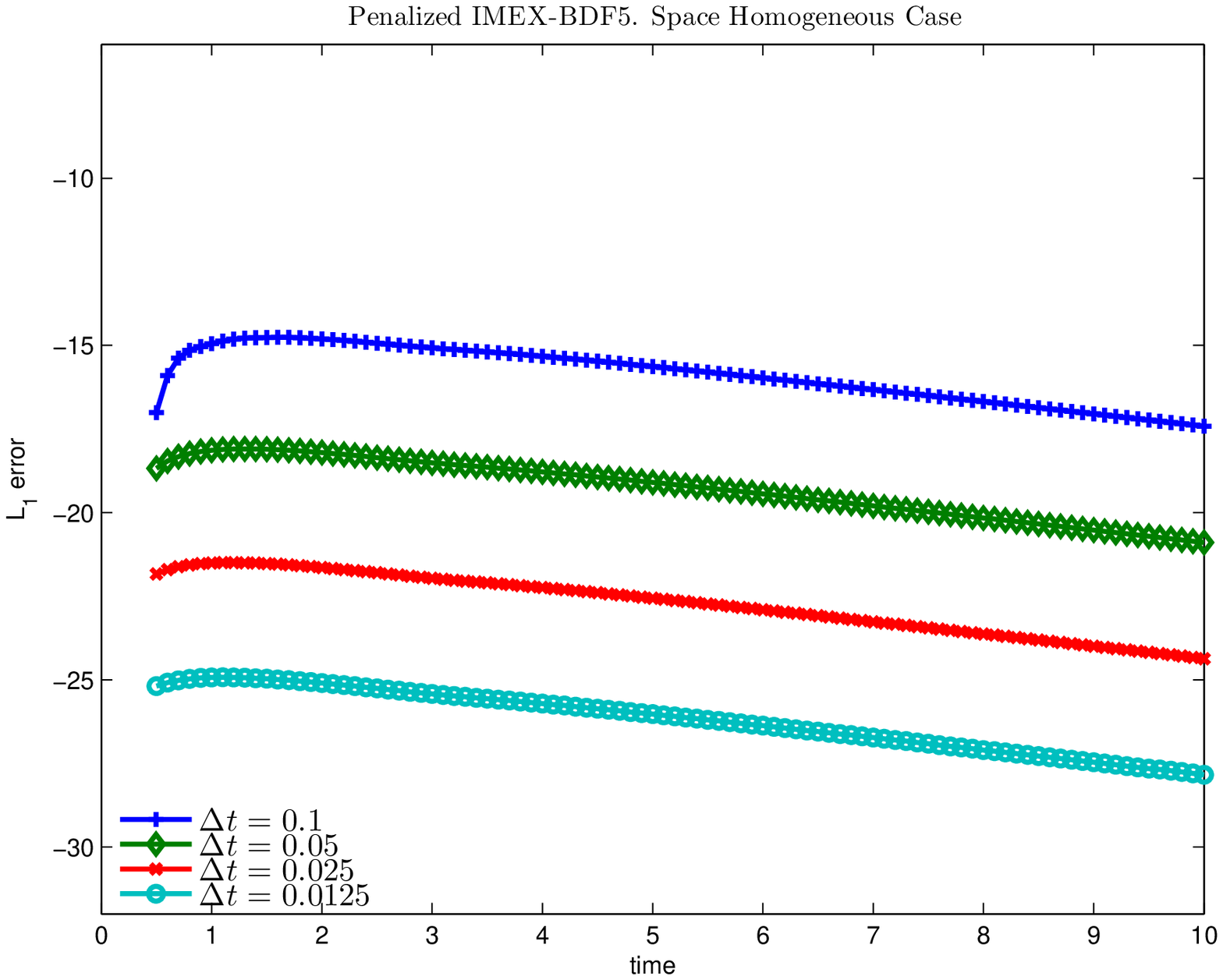}
\caption{$L_{1}$ error for the distribution function $f$. Penalized IMEX-BDF schemes. From left to right, top to bottom: second, third, fourth and fifth order methods.}\label{fig:conv_Bolt_H}
\end{center}
\end{figure}

\subsection{Homogeneous Boltzmann equation}
This test has been performed to study the accuracy of the penalized IMEX multistep schemes based on the BGK operator in the case of the space homogeneous Boltzmann equation
\be
\partial_t f 
=\frac1{\varepsilon}Q(f,f),
\ee
with $Q(f,f)$ given by (\ref{eq:Q}) in the two dimensional velocity space. The collision kernel corresponds to Maxwellian molecules $B(|v-v_*|,n)=4\pi$ and the
fast spectral method \cite{MP} has been used to approximate the collision operator with
$N_{v}=64$ grid points in each velocity direction and a grid
$[-v_{\max}, v_{\max}]^2$ with $v_{\max}=10$. 

As discussed in Section \ref{sec:stab} the choice of the penalization factor plays a major rule in the case of high order methods. Here we assume a constant $\mu=\rho$ so that the loss part of the collision term cancel in the penalization process and we have
\[
Q(f,f)-\mu(M-f)=Q^+(f,f)-\rho f - \mu (M-f) =Q^+(f,f)-\mu M,
\]
where $Q^+$ is the gain part of the collision operator (\ref{eq:Q}).
 
The non equilibrium
initial data is given by 
\[
f(v,0)=\frac{\rho_0}{4\pi T_0}\left(\exp\left(\frac{(v_x-u_x)^2+(v_y-u_y)^2}{2T_0}\right)+\exp\left(\frac{(v_x+3u_x)^2+(v_y-u_y)^2}{2T_0}\right)\right),
\]
where we took $\rho_0=1$, $u_x=1$, $u_y=1$ and $T_0=1$. 
We report the results in term of the $L_1$ error for the distribution function $f$, i.e. $\|f-f_{ex}\|_1$, for the penalized IMEX-BDF schemes in Figure \ref{fig:conv_Bolt_H}. Analogous results are obtained using the other IMEX multistep schemes of the same order.   
The reference solution has been computed by using a third order penalized Runge-Kutta method \cite{dimarco7} using very small time steps up to machine error. The same third order penalized IMEX Runge-Kutta scheme has been employed also
to produce the initial vector needed for starting the multistep schemes. The same strategy has also been adopted: the time step has been reduced up to the machine error for producing numerically exact initial values.
The theoretical convergence rate have been observed for all the schemes considered for this test in all regimes. In Figure \ref{fig:conv_Bolt_H}) we report the error curves for different values of the time step for the penalized IMEX-BDF methods. Essentially the same error curves are obtained with the other IMEX multistep schemes of corresponding order. Note that in this case we have only one time scale $t/\varepsilon$ and therefore we fixed $\varepsilon=1$ and consider different values of the time step. A similar test case was considered in \cite{dimarco8} for penalized IMEX Runge-Kutta methods.

\subsection{Non homogeneous Boltzmann equation}
The last test case considers the numerical solution of the full non homogeneous Boltzmann equation (\ref{eq:1b})-(\ref{eq:Q})
in one space dimension and in the two dimensional velocity space. The same penalized setting of the homogeneous test has been used, namely Maxwell molecules as a collision kernel and BGK model as penalization operator with $\mu(x,t)=\rho(x,t)$. 

\begin{figure}
\begin{center}
\includegraphics[scale=0.35]{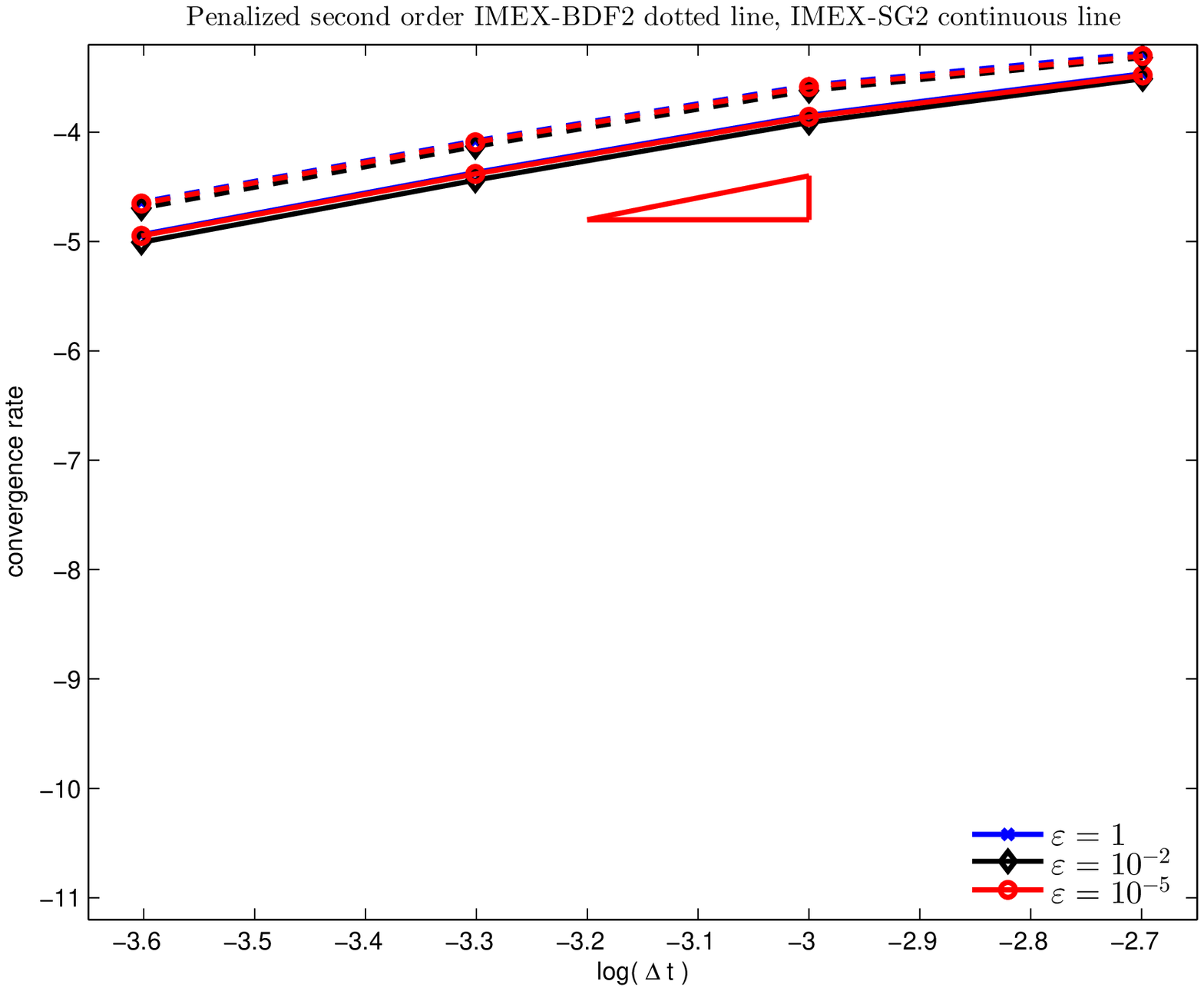}
\includegraphics[scale=0.35]{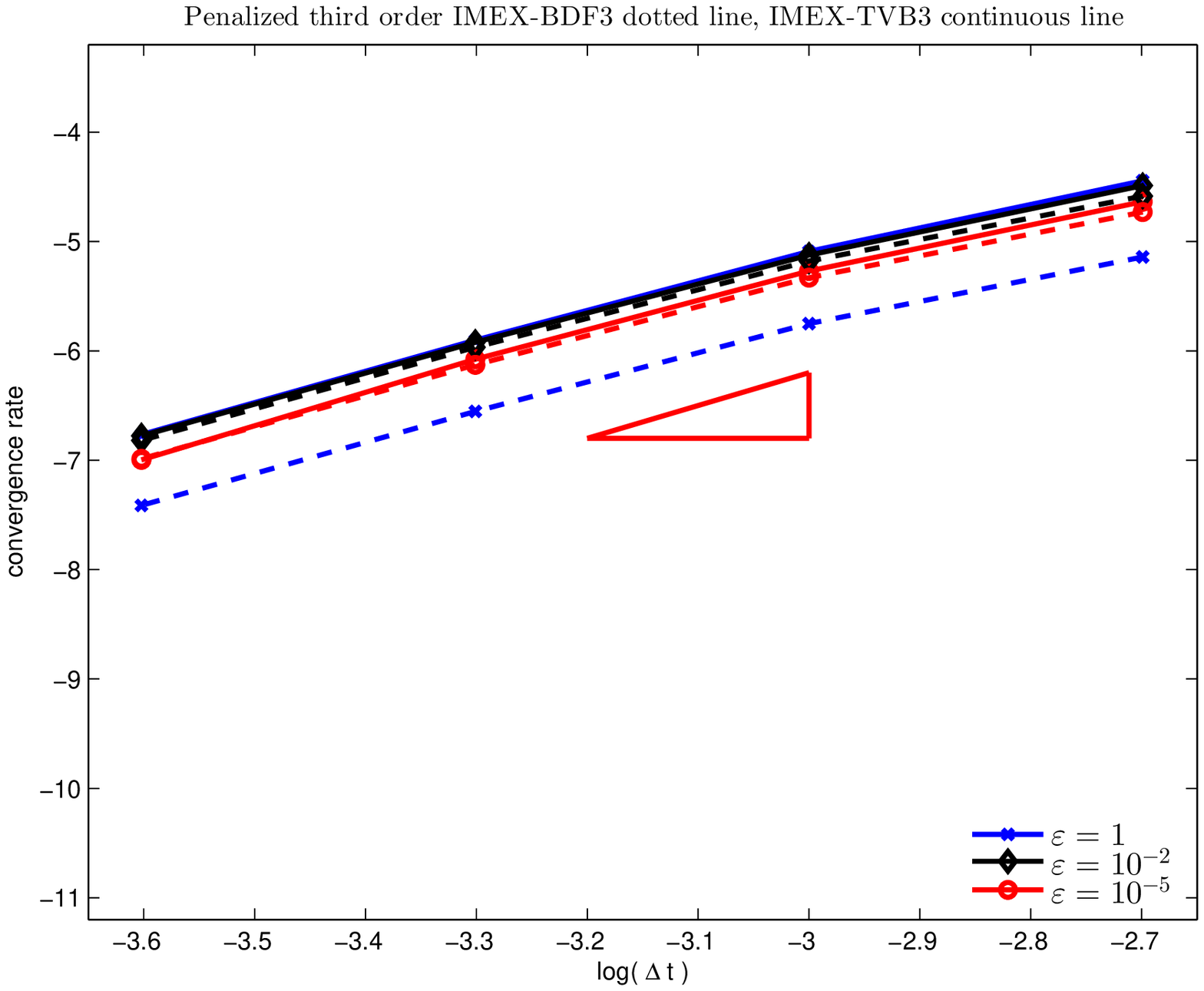}\\
\includegraphics[scale=0.35]{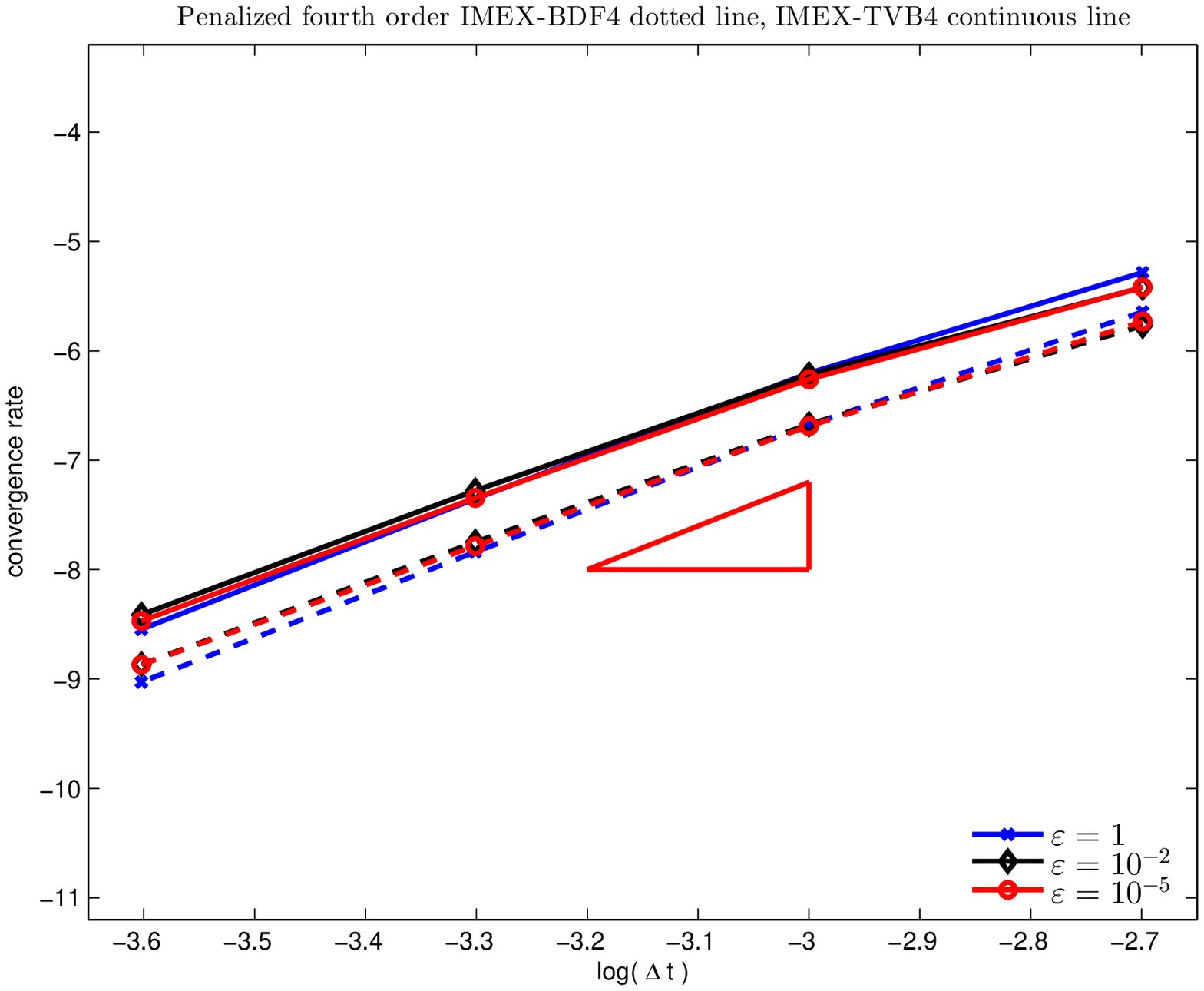}
\includegraphics[scale=0.35]{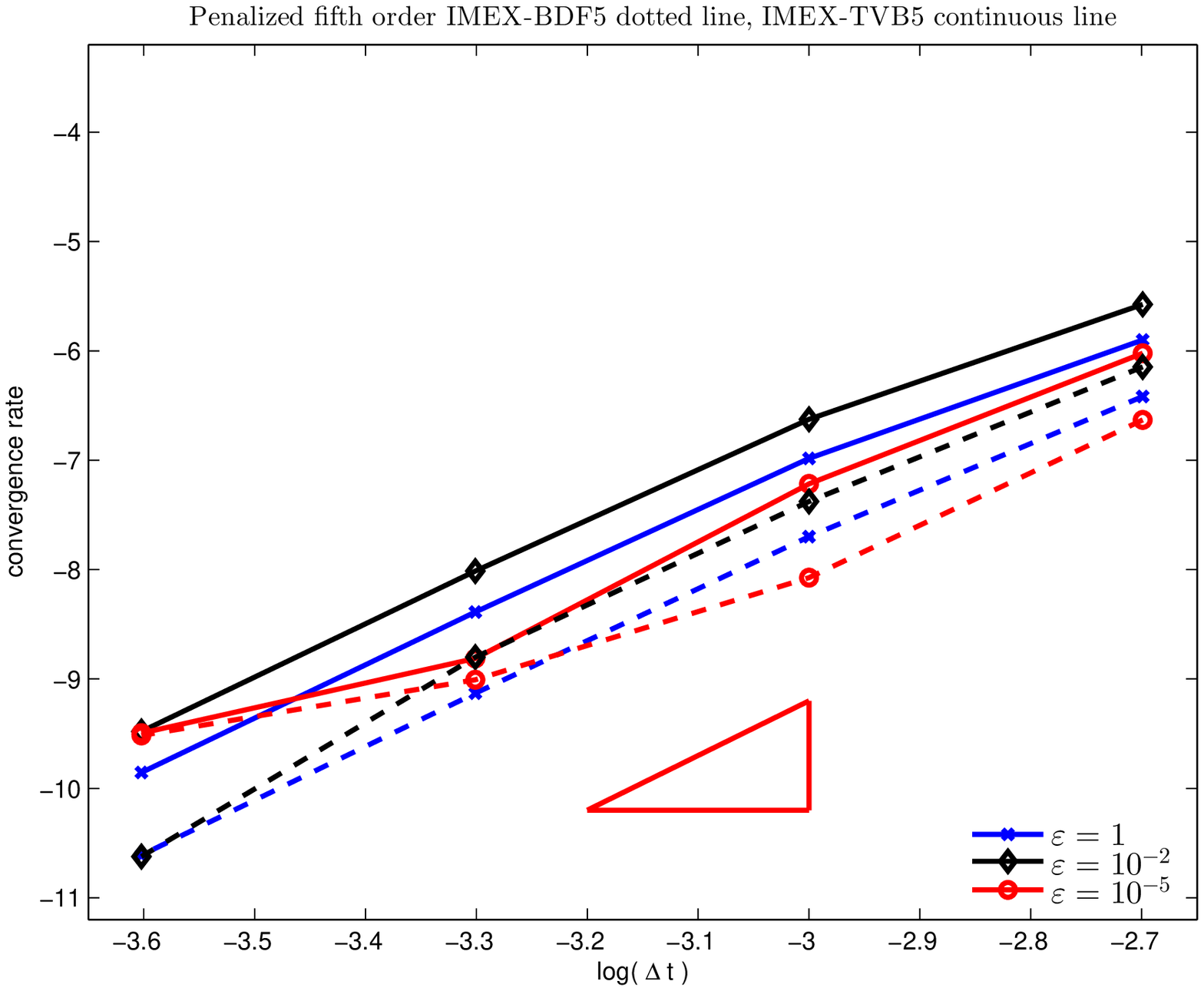}
\caption{Convergence rate for the density $\rho$ for the penalized IMEX-BDF2 and IMEX-SG2 schemes (top left), penalized IMEX-BDF3 and IMEX-TVB3 schemes (top right), 
penalized IMEX-BDF4 and IMEX-TVB4 schemes (bottom left) and penalized IMEX-BDF5 and IMEX-TVB5 schemes (bottom right).}\label{fig:conv_Bolt1}
\end{center}
\end{figure}

The computation is performed on $(x, v) \in [0, 1] \times
[-v_{\max}, v_{\max}]^2$, with $v_{\max} = 8$ and $N_v = 256$. A $5$-th order WENO
scheme has been used for the space discretization with the largest time step equal to $\Delta t_{max} = \frac{\Delta x}{4v_{\max}}$ and the initial data  \be
\varrho_{0}(x)=\frac{2+\sin(8\pi x)}{3}, \ u_{0}(x)=0, \ T_{0}(x)=\frac{2+\cos(8 \pi x)}{3}.\ee 
As in the case of the space non homogeneous BGK model an initial distribution $f(x,v,t=0)=f_0=M[f_0]+\varepsilon g$ with $g$ a perturbation from equilibrium consistent with the Navier-Stokes limit has been used.
We report the convergence rates for the multistep schemes obtained by measuring the $L_1$ norm of the error for the density for different values of the
Knudsen number, i.e. $\varepsilon=10^{-1}$, $\varepsilon=10^{-2}$ and $\varepsilon=10^{-5}$ and different time steps. 
The rates, as before, have been obtained repeating the computation with the same initial data for a decreasing time step
$\Delta t_{1} = \Delta t_{max}/2$, $\Delta t_{2} = \Delta t_{max}/4$  and $\Delta t_{3} = \Delta t_{max}/8$
while the grid points in space remain fixed to $N_x=128$. The penalized multistep schemes are initialized by a penalized third order IMEX Runge-Kutta scheme reaching the machine precision.
Figure \ref{fig:conv_Bolt1} show the convergence rates for various penalized IMEX multistep schemes from order two up to order five. The schemes tested exhibit the theoretical convergence rate, except for the fifth order schemes where degradation of accuracy is observed for small values of $\varepsilon$. This is particularly evident for IMEX-BDF5 scheme and it is due to the fact that very high order penalized IMEX multistep schemes are extremely sensitive to the choice of the penalization factor, as analyzed in Section \ref{sec:stab}. 
Therefore, one should consider a better penalization strategy, for example based on a velocity dependent frequency BGK model. An interesting alternative in this case is represented by the ES-BGK model \cite{BPESBGK}. Here, however, we do not explore further this direction.

\section{Final considerations}
IMEX multistep schemes represent an interesting alternative to IMEX Runge-Kutta and exponential methods, in particular when one deals with multidimensional kinetic equations with stiff collision terms. Thanks to their simpler structure, it is possible to achieve high order accuracy without additional coupling conditions and at a reduced computational cost. This is of paramount importance in the numerical solution of kinetic equations, where the computational cost due to the multidimensionality of the problem and the structure of the Boltzmann collision term is often the major concern in the construction of the numerical method. In this paper we have analyzed and studied such schemes in the case of general collisional kinetic equations, by considering a penalized strategy to avoid the costly inversion of the collision term in the case of the full Boltzmann equation. Both, the theoretical results and the numerical tests have shown that the schemes are able to work with uniform accuracy for a wide range of the relaxation parameter $\varepsilon$ including the Navier-Stokes regime. On the other hand the schemes are more sensitive to the choice of the penalization factor and instabilities may be observed with high order penalized methods. Further research directions will consider the possibility to use high order schemes without penalization by a suitable inversion of the single collision term needed at each time step and the use of velocity dependent frequencies in the BGK model used for penalization, like the ES-BGK model.





\begin{thebibliography}{99}
\small

\bibitem{Ascher} {\sc U.~M.~Ascher, S.~J.~Ruuth, R.~J.~Spiteri}, {\it Implicit-explicit Runge-Kutta methods for time dependent
Partial Differential Equations}, Appl. Numer. Math., 25 (1997),
151--167.

\bibitem{Ascher2} {\sc U.~M.~Ascher, S.~J.~Ruuth, B.~T.~R.~Wetton}, {\it Implicit-Explicit methods for time-dependent partial differential equations}, SIAM J. Num. Anal. 32, (1995), 797--823.

\bibitem{BLM} {\sc M.~Bennoune, M.~Lemou, L.~Mieussens}, {\it Uniformly stable numerical schemes for the Boltzmann equation
preserving the compressible Navier-Stokes asymptotics}, J. Comp.
Phys., 227 (2008), 3781--3803.

\bibitem{BGK}
{\sc P. L. Bhatnagar, E. P. Gross, M. Krook,} {\it A model for
collision processes in gases I: Small amplitude processes in charged
and neutral onecomponent systems}, Phys. Rev. 94, 511 (1954).

\bibitem{BJ}
{\sc B. Yan, S. Jin,} {\it A successive penalty-based asymptotic-preserving scheme for kinetic equations}, SIAM J. Sci. Comp., 35, 150-172, (2013).





\bibitem{boscarino}{\sc S.~Boscarino, G.~Russo}, {\it
On a class of uniformly accurate IMEX Runge-Kutta schemes and
applications to hyperbolic systems with relaxation},  SIAM J. Sci.
Comput., 31 (2009), 1926--1945.

\bibitem{BGP00} 
{\sc F.~Bouchut, F.~Golse, M.~Pulvirenti}, {Kinetic equations and asymptotic theory}, {Gauthiers-Villars} (2000).

\bibitem{BPESBGK} {\sc F.~Bouchut, B.~Perthame}, {\it A BGK model for small Prandtl numbers in the Navier-Stokes approximation}, J. Stat. Phys. 71, (1993), 191--207.


\bibitem{Caflisch} 
{\sc R.~E.~Caflisch}, {\it The fluid dynamic limit of the nonlinear Boltzmann equation}, Commun.
Pure Appl. Math. 33 (1980), 651--666.

\bibitem{CJR}
{\sc R.~E.~Caflisch, S.~Jin, G.~Russo}, {\it Uniformly accurate
schemes for hyperbolic systems with relaxation}, SIAM J. Numer.
Anal., 34 (1997), 246--281.




\bibitem{cercignani}
{\sc C.~Cercignani}, {The Boltzmann equation and its
applications}, Springer-Verlag, New York, (1988).





\bibitem{degond1}
{\sc P.~Degond, S.~Jin, L.~Mieussens}, {\it A smooth transition
between kinetic and hydrodynamic equations}, J. Comp. Phys., 209
(2005), 665--694.




\bibitem{dimarco5}
{\sc P. Degond, G. Dimarco, L. Mieussens}, {\it A multiscale
kinetic-fluid solver with dynamic localization of kinetic effects.}
J. Comp. Phys., 229 (2010), pp. 4907-4933.


\bibitem{dimarco2}
{\sc G.~Dimarco, L.~Pareschi}, {\it Fluid solver independent
hybrid methods for multiscale kinetic equations}, SIAM J. Sci.
Comput., 32 (2010), 603--634.

\bibitem{dimarco6}
{\sc G.~Dimarco, L.~Pareschi}, {\it Exponential Runge-Kutta methods
for stiff kinetic equations.} SIAM J. Num. Analysis, 49  (2011), pp.
2057--2077.

\bibitem{dimarco7}
{\sc G.~Dimarco, L.~Pareschi}, {\it Asymptotic-preserving IMEX Runge-Kutta methods for nonlinear kinetic equations}, SIAM J. Num. Anal. (2013), 1064--1087.

\bibitem{dimarco8}
{\sc G. Dimarco, L. Pareschi}, {\it High order asymptotic preserving
schemes for the Boltzmann equation}, Comptes Rendus Mathematique 350, 9, (2012), 481--486. 

\bibitem{Acta}
{\sc G.~Dimarco, L.~Pareschi}, {\it Numerical methods for kinetic equations}, Acta Numerica 23, (2014) 369--520.

\bibitem{DB}
{\sc D.~R.~Durran, P.~N.~Blossey}, {\it Implicit-Explicit Multistep Methods for Fast-WaveÐSlow-Wave Problems}, Mon. Wea. Rev. 140, (2012), 1307--1325.



\bibitem{Filbet} {\sc F.~Filbet, S.~Jin}, {\it A class of asymptotic preserving schemes for kinetic equations and related
problems with stiff sources}, J. Comp. Phys., 229 (2010), pp.
7625-7648.


\bibitem{Filbet3} {\sc F.~Filbet, S.~Jin}, {\it An asymptotic preserving scheme for the ES-BGK model of the Boltzmann equation}.
J. Sci. Comput., 46 (2011), 204-224.

\bibitem{FHV}
{\sc J.~Frank, W.~Hundsdorfer, J.~G.~Verwer}, {\it On the stability of implicit-explicit linear multistep methods}. Appl. Num. Math. 25, (1997), 193--205.

\bibitem{toscani}
{\sc E.~Gabetta, L.~Pareschi, G.~Toscani}, {\it Relaxation schemes
for nonlinear kinetic equations}, SIAM J. Numer. Anal., 34 (1997),
2168--2194.

\bibitem{Golse}
{\sc F.~Golse}, {\it The Boltzmann equation and its hydrodynamic limits}, in {Handbook of Differential equations}, {Evolutionary Equations}, vol. 2. Edited by C.M.~Dafermos and E.~Feireisl, Elsevier (2005).

\bibitem{Gosse}
{\sc L.~Gosse}, {Computing Qualitatively Correct Approximations of
  Balance Laws. {E}xponential-Fit, Well-Balanced and Asymptotic-Preserving}, SEMA SIMAI Springer Series (2013).


\bibitem{GRS} 
{\sc M.~Groppi, G.~Russo, G.~Straquadanio}, {\it High order semi-Lagrangian methods for the BGK equation}, Comm. Math. Sci. (to appear).

\bibitem{HairerWanner}
{\sc E.~Hairer, G.~Wanner}, {Solving Ordinary Differential
Equations {\rm II:}} {Stiff and Differential-Algebraic Problems},
Springer-Verlag, New York, (1987).







\bibitem{Hol} {\sc L.~H.~Holway}, {\it Kinetic theory of shock structure using an ellipsoidal
distribution function}, Rarefied Gas Dynamics, I, Academic Press
(1966), 193--215.

\bibitem{HLP} {\sc J. Hu, Q. Li, L. Pareschi}, {\it Asymptotic-preserving exponential methods for the quantum Boltzmann equation with high-order accuracy}, J. Sci. Comput., 62 (2015), pp. 555-574



\bibitem{HR}
{\sc W.~Hundsdorfer, S.~J.~Ruuth}, {\it
IMEX extensions of linear multistep methods with
general monotonicity and boundedness properties}, J. Comp. Phys. 225, (2007), 2016--2042.

\bibitem{HRS}
{\sc W. Hundsdorfer, S.J. Ruuth, R.J. Spiteri}, {\it Monotonicity-preserving linear multistep methods}, SIAM J. Numer. Anal. 41, (2003) 605--623.


\bibitem{Jin2}
{\sc S.~Jin}, {\it Efficient Asymptotic-Preserving (AP) schemes for some multiscale kinetic equations},
SIAM J. Sci. Comput., 21 (1999), 441--454.


\bibitem{Lem}
{\sc M.~Lemou},  {\it Relaxed micro-macro schemes for kinetic
equations}, C. R. Acad. Sci. Paris, Ser. I, 348 (2010), 455--460.

\bibitem{Qin}
{\sc Q.~Li, L.~Pareschi},  {\it Exponential Runge-Kutta schemes for inhomogeneous Boltzmann equations with high order of accuracy}. J. Comp. Phys. 259 (2014) 402--420. 



\bibitem{MP}
{\sc C.~Mouhot, L.~Pareschi}, {\em Fast algorithms for computing the
Boltzmann collision operator}, Math. Comp. 75 (2006),
pp.~1833--1852.



%

\bibitem{PRimex}
{\sc L.~Pareschi, G.~Russo}, {\it Implicit-Explicit Runge-Kutta
methods and applications to hyperbolic systems with relaxation},
J. Sci. Comput., 25 (2005), 129--155.
%
\bibitem{Puppo}
{\sc S.~Pieraccini, G.~Puppo}, {\it Implicit-Explicit schemes for
BGK kinetic equations}, J. Sci. Comput., 32 (2007), 1-28.





\bibitem{Shu} {\sc C.-W. Shu}, \emph{High Order Weighted Essentially Nonoscillatory Schemes for Convection Dominated Problems} 
SIAM Rev. 51, pp. 82-126 (2009).

\bibitem{KT}
{\sc S. Tiwari, A. Klar}, {\it An adaptive domain decomposition procedure for Boltzmann and Euler equations}, J. Comp. and Appl. Math., 90, (1998), 223--237.



\end{thebibliography}
\end{document}